\documentclass[12pt]{article}
\usepackage{makeidx}
\usepackage{epsfig}
\usepackage{float}
\usepackage{amscd}
\usepackage{amsmath}
\usepackage{amssymb}
\usepackage{color}
\usepackage{fancyheadings}
\usepackage{fancyvrb}
\usepackage[]{algorithm2e}
\RestyleAlgo{ruled}

\newcommand{\intox}[1]{\int_{\Omega} #1 \, dx}
\newcommand{\nn}{{n}}
\newcommand{\bn}{\mathbf{\nn}}
\newcommand{\uu}{{u}}
\newcommand{\du}{{Du}}
\newcommand{\vv}{{v}}
\newcommand{\ww}{{w}}
\newcommand{\xx}{{x}}

\newcommand{\sdiv}{{{\nabla\cdot} \,}}
\newcommand{\bfz}{{\bf 0}}
\newcommand{\gbc}{{\bf g}}

\newcommand{\omitit}[1]{}
\newcommand\half{{\textstyle{\frac{1}{2}}}}
\newcommand{\set}[2]{\left\lbrace #1 \; : \; #2 \right\rbrace}
\newcommand{\Reyuls}{{\Bbb R}}

\date{}

\makeindex
\usepackage{makeidx}

\newtheorem{theorem}{Theorem}[section]
\newtheorem{remark}[theorem]{Remark}

\begin{document}
\title{Chaotic dynamics of two-dimensional flows \\ around a cylinder}
\author{
L. Ridgway Scott\\ Emeritus, University of Chicago \\
Rebecca Durst \\ Department of Mathematics, University of Pittsburgh }
 \def\thepage{}\maketitle\pagenumbering{arabic}
\centerline{\today}

\begin{abstract}
We study flow around a cylinder from a dynamics perspective, using drag and
lift as indicators.
We observe that the mean drag coefficient bifurcates from the steady case
when the Karman vortex street emerges.
We also find a jump in the dimension of the drag/lift attractor just above
Reynolds number 100.
We compare the simulated drag values with experimental data obtained
over the last hundred years.
Our simulations suggest that a vibrational resonance in the cylinder would
be unlikely for Reynolds numbers greater than 1000, where the drag/lift
behavior is fully chaotic.
\end{abstract}

\section{Introduction}

Flow around a cylinder has long been a test problem of interest in fluid dynamics
\cite{ref:dryden1930pressure,ref:dragecrisisfage}.
This problem has been examined extensively both experimentally
\cite[Figure 14.15]{PantonRonaldL2013IF}
and computationally (see references in \cite{lrsBIBiw,lrsBIBjn}).
It is well known that a Hopf bifurcation occurs near Reynolds number 50
\cite{ref:cylinderHopfvortexstreet,chen1995bifurcation},
where steady flow gives way to the
Karman vortex street \cite{birkhoff1953formation,ref:karmanstreetexpts}.
Recent computational studies \cite{ref:Rajani2D3Dcylinderflow,ref:confinedcylinderflow} 
have provided detailed information of significant interest.  Here we extend such computational studies with a focus on the dynamic behavior of the flow for Reynolds numbers up to $\approx 10^4$.

We limit our study to two-dimensional flow to keep the length manageable.
It has been known for some time, both computationally
\cite{ref:BeaudanMointhesis94,ref:dongkarni05dns}
and experimentally \cite{ref:FazleHayakawatreedee}, that flow past a cylinder
remains largely two-dimensional up to a Reynolds number of at least $10^4$.
Indeed, \cite{ref:FazleHayakawatreedee} studies the three-dimensional
deviation from two-dimensional flow.
Many issues of interest are three-dimensional, but it is useful to
see exactly what features can be approximated as two dimensional.
Further, the two dimensional flow problem is a well posed mathematical problem
of independent interest.
Applications of two-dimensional flow are discussed in \cite{ref:wolfchaocylinderflow}.

Our goal for this paper is to provide computational insight into the dependence of the flow dynamics \cite{wiggins2003introduction} on the Reynolds number. We quantify these dynamics with standard metrics: the Lyapunov exponent, the Strouhal number, and the fractal dimension. Additionally, we provide references to significant physical experiments to serve as a comparison to our computational experiments. These are detailed in Section \ref{sec:expData}.

Our computational technique involves so-called pressure robust
finite elements \cite{ref:LinkeBeltramiFlows} for the spatial discretization
and IMEX (linearly IMplicit, nonlinearly EXplicit) time stepping \cite{lrsBIBis},
as explained in section \ref{sec:imexeme}.
This utilizes a spatial discretization that is essentially divergence-free.
The time stepping requires no stabilization, provided the divergence of the
approximate solution is kept small enough \cite{lrsBIBis} and the time step
is also small enough.

 One result of this paper is that we observe strong evidence from multiple metrics that the vortex shedding in the Karm\'an vortext street is periodic, beginning around Reynolds number $50$ (coinciding with the known Hopf bifurcation \cite{ref:cylinderHopfvortexstreet,chen1995bifurcation}) and continuing this periodicity up to Reynolds number $200$.  Moreover, our results indicate that this periodicity noticeably begins to break down as early as Reynolds number $250$.  We know of no mathematical proof of the existence of a periodic solution at the lower Reynolds numbers.
 
 Additionally,  we observe a bifurcation in the mean drag value as the flow goes from steady to time-dependent near Reynolds number 50, where the Karm\'an vortex street emerges.  {We expect a change in the flow solution from steady to unsteady at the point, however we note a bifurcation in the \textit{mean drag} indicates that the drag value of the time-dependent flow does not oscillate around the drag value of the steady-state solution at the same Reynolds number. In other words, the steady state solution is not the average of the time-dependent solution.}

\subsection{Experimental data}\label{sec:expData}

Experimental data for the drag on a cylinder has been reported for over
a century \cite{ref:errelflowRex,wieselsberger1921neuere}. For example,  in \cite{ref:midRecylinderdrag}, data from several papers are collated, including
data from \cite{ref:errelflowRex}.
References \cite{ref:DelaneySorensenNACAcyldrag,ref:cylendardrag} report
on  Reynolds numbers $R\in[10^4,10^7]$.
In particular, \cite[Figure 1]{ref:cylendardrag} reproduces
\cite[Figure 5]{ref:DelaneySorensenNACAcyldrag}.
These figures indicate how cylinders of different sizes are used to determine
drag coefficients in different ranges of Reynolds numbers.
Comparing the results for cylinders of size 4 and 1 (inches) in their
overlapping range of Reynolds numbers ($3\times 10^4$ to $10^5$) gives
a sense of the experimental error.

There are two physical experiments that serve as valuable reference data to compare the results and validity of our numerical scheme. The first is data collected by Relf in 1914 \cite{ref:errelflowRex}, reporting the drag force measured on wires of small diameter. This data from Relf is reproduced in
\cite[Figure 6]{ref:midRecylinderdrag} for Reynolds numbers $R\in[1,10^3]$.
That figure compares the Relf data with data from \cite{wieselsberger1921neuere,
ref:WieselsbergerNACA}.
For Reynolds numbers $R\in[10^3,3\times 10^4]$,
the data from Relf \cite{ref:errelflowRex} is reproduced in
\cite[Figure 7]{ref:midRecylinderdrag}.
The data in \cite[Figure 1]{ref:WieselsbergerNACA} is also represented
in these two figures in \cite{ref:midRecylinderdrag}.
Reference \cite{ref:WieselsbergerNACA} is derived from
\cite{wieselsberger1921neuere}.  Relf's approach was to measure the force on (1) a frame of very fine wires
for small Reynolds numbers and (2) solitary (bigger) wires
for larger Reynolds numbers.
In Figure \ref{fig:relfdata} below, we have plotted the different groups
of measurements for different wire sizes in different colors, with the data
points connected by straight lines.
This also provides a useful guide to the size of the experimental
uncertainty that we expect to be present in experiments conducted more than a century ago. However, despite concerns over the age of the physical data, we still found it to be a remarkably reliable reference for our numerical method.

The second source of physical data comes from experiments published by Tritton in 1959 \cite{ref:trittonlowrexpts} 
for low Reynolds numbers by measuring flow past quartz crystal cylinders. In \cite[page 553]{ref:trittonlowrexpts}, the
experimental error in the drag coefficient is estimated to be around 6\%,
although the data scatter suggests it may be closer to 10\%.
It is also noted in \cite[page 554]{ref:trittonlowrexpts} that cylinder vibration
may affect drag observations.
The comparison between experiments and simulation data is within the experimental error,
except for $300\leq R\leq 4000$ or so.
The discrepancy in this range could be due to vibrations of the
cylinders (wires) in the experiments \cite{ref:annurevfluidvortexvibrations}.
In \cite[Figure 6.9, page 179]{ref:Meneghini1994Thesis}, it is shown that 
forced vibrations affect the drag coefficient substantially.

\subsection{Outline of the paper}
The remainder of this paper proceeds as follows. In Section \ref{timeDNS} we introduce the time-dependent Navier-Stokes equations and the domain for our experiments. Additionally, in this section we provide a description of our computational method.  Subsequently, in Section \ref{dragLift} we introduce the drag and lift coefficients and provide an initial comparison of our numerical computations with physical experiements. Finally, in Section \ref{dyn} we provide a dynamical analysis of our drag and lift computations using a variety of dynamical metrics.

\section{Time-dependent Navier-Stokes} \label{timeDNS}

For incompressible fluid flow,  the Navier-Stokes equations are as follows,
\begin{equation} \label{eqn:nastoensn} 
\uu_t-\nu\Delta\uu+\uu\cdot\nabla\uu+\nabla p = 0,\qquad \sdiv\uu=0,\quad\hbox{in}\;\Omega,
\end{equation}
 where $\nu$ is the kinematic viscosity.
 
We define the initial condition $\uu(t=0)=\uu_0$,
and set the boundary condition $\uu=\gbc$ on $\partial\Omega$.

Let $V$ be the subset of the $H^1(\Omega)^2$ consisting of divergence-free functions.  Here, $H^1(\Omega)^2$ is the Sobolev space consisting of square integrable functions whose gradients are square integrable. Then we can solve \eqref{eqn:nastoensn} via the variational formulation
\begin{equation} \label{eqn:bisneqvst}
\intox{\vv_t\cdot\ww} +{\nu}\intox{\nabla\vv:\nabla\ww} +  c(\vv,\vv,\ww)=0
\end{equation}
for all $\ww\in V$, where $c(\uu,\vv,\ww)$ is defined by
\begin{equation} \label{eqn:seevarformqvst} 
c(\uu,\vv,\ww)=\intox{(\uu\cdot\nabla\vv)\cdot\ww}.
\end{equation}

\subsection{Computational domain and boundary conditions}

\begin{figure}[H]
\centering
\includegraphics[height=0.3\textwidth]{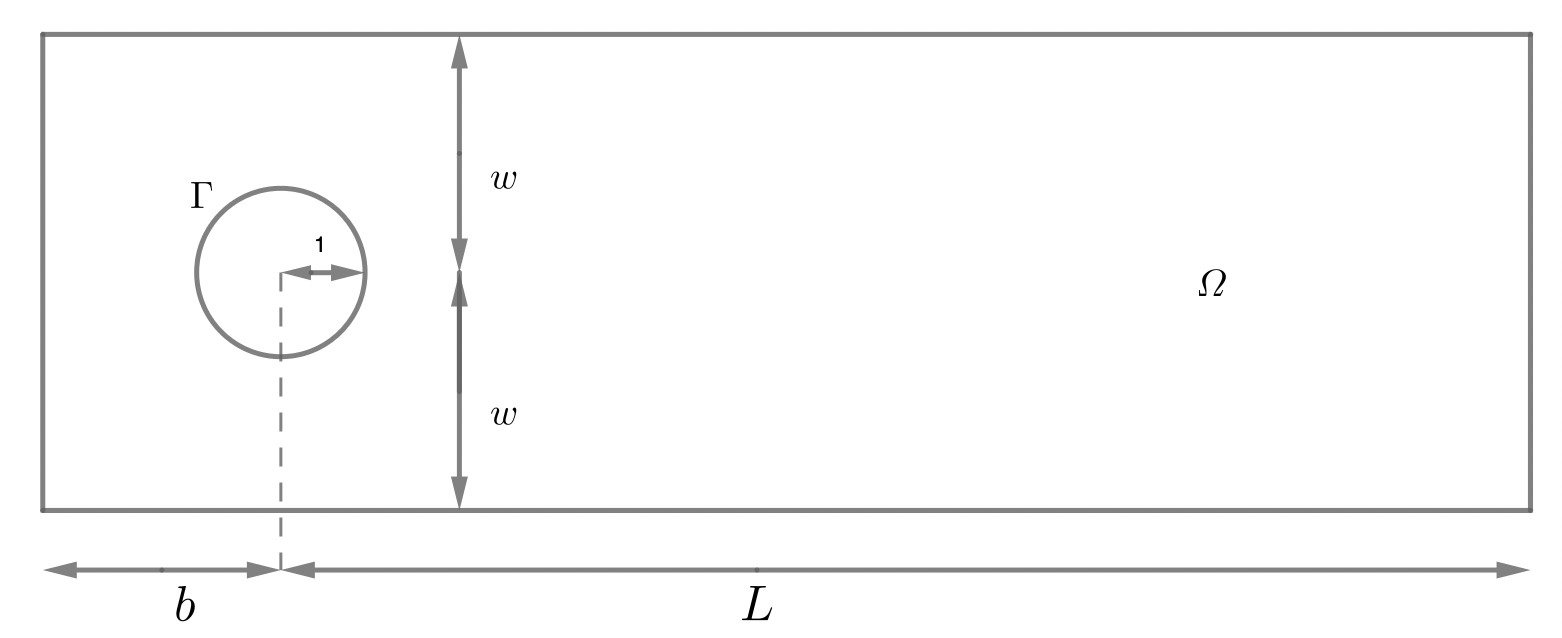}
\caption{The domain $\Omega$.}\label{omeg}
\end{figure}

The computational domain $\Omega$ (see Figure \ref{omeg}) is defined by various parameters:
\begin{equation} \label{eqn:domparam} 
\Omega=\set{(x,y)}{-b<x<L,\;|y|<w,\; x^2 + y^2 >1}.
\end{equation}
This computational domain surrounds a cylinder with radius $1$, centered at $(0,0)$. The boundary of this cylinder is $\Gamma$, defined by:
$$
\Gamma=\set{(x,y)}{x^2+y^2=1}.
$$

The Reynolds number is typically defined using the cylinder diameter as the length
scale, so $R=2/\nu$ if we assume that
the maximum inflow speed is 1, as we do in all computations here.

The choice of $b$, $L$, and $w$ are somewhat arbitrary, and so we will refer to the
rectangle that is the boundary of $[-b,L]\times[-w,w]$ as the {\em computational boundary}.  However,  we note that in \cite{lrsBIBiw},  extensive computations were done with $b=w=6$, $L=12$, and $d=1$.
These suggest that Reynolds--Orr time-dependent instabilities appear for 
Reynolds numbers significantly below 20, and thus well below the threshold for the
emergence of the periodic Karman vortex street,
near $R=50$.

We consider \textit{free-stream boundary conditions}: $\gbc=(1,0)$,
 on the computational rectangle,  but with $\gbc=\bfz$ on $\Gamma$.

\subsection {Numerical methods} \label{sec:imexeme}
Let $W$ be a subset of $H^1(\Omega)^2$. In this paper,  we model this system using a second-order,  linearly \textit{implicit},  nonlinearly \textit{explicit} (IMEX) time-stepping scheme \cite[section 8.6]{lrsBIBis} given by
the variational formulation
\begin{equation} \label{eqn:secordhersiml}
\begin{split}
\frac{1}{\Delta t}\intox{(\uu^{n+1}-\vv^n)\cdot\ww} &+ \frac{\nu}{2}
\intox{\big(\nabla\uu^{n+1}+\nabla\uu^{n}\big):\nabla\ww}
+ c(E\uu^{n},E\uu^{n},\ww)=0
\end{split}
\end{equation}
for all $\ww\in W$, where $E\uu^{n}$ is obtained by extrapolation:
$$
E\uu^{n}=\frac32 \uu^n - \frac12 \uu^{n-1}\approx \uu^{n+\half},
$$
where we define $\vv^0$ to be equal to the stationary Stokes solution,  and set $E\vv^0 = \vv^0$. 

In practice,  this method is implemented using Scott-Vogelius elements of degree 4, so we let $W$ be the subset of $H^1(\Omega)^2$ consisting of piecewise quartic polynomials.  Issues at the boundary can occur due to the polygonal approximation of the cylinder,  therefore we enforce the boundary conditions on the cylinder weakly using Nitsche's method,  as explained in \cite{lrsBIBkj}.  Additionally,  the incompressibility constraint is enforced with the iterated-penalty method \cite{lrsBIBih} with $\rho$ proportional to the time-step $\tau$.

More specifically, let $n\geq1$ and assume the solutions at $(\uu^n, p^n)$ are known.  For  $\tau=2/(\nu\Delta t)$, we define a parametrized variational form 
\begin{alignat}{1}
a^\tau(\vv,\ww) =& \tau\int_{\Omega}\vv \cdot \ww \, dx  + \int_{\Omega} \nabla \vv : \nabla \ww \, dx + \rho \int_{\Omega} (\sdiv{\vv})(\sdiv{\ww}) \, dx \nonumber \\
& + \gamma \int_{\Gamma} \vv \cdot \ww \, ds - \int_{\Gamma} (\nabla \vv \,\bn) \cdot \ww \, ds- \int_{\Gamma} \vv \cdot (\nabla \ww\, \bn ) \, ds, \label{aTau}
\end{alignat}
and a linear form depending on the solution $\uu^n$,
\begin{alignat}{1}\label{fDef}
F^{n}(\ww)&=\tau\,(\uu^{n},\ww)-a^0(\uu^{n},\ww)-\frac{2}{\nu}c(E\uu^{n},E\uu^{n},\ww).
\end{alignat}

 We note that the parameter $\gamma$ above is the Nitsche penalization parameter, and $\bn$ is the outward facing unit normal to the fluid domain, $\Omega$ along the boundary $\Gamma$. Consequently,  we solve for the updated solution $(\uu^{n+1}, p^{n+1})$ with the following algorithm:

\begin{algorithm}[H]
\caption{Scott-Vogelius-Nitsche method with IMEX time-stepping}\label{alg:svnIMEX} 
\begin{enumerate}
\item[(1)] Set $\vv_0 = \mathbf{0}$ and $\epsilon>0$, and run the following iterative method:

\While{$\|\sdiv{\uu_k}\|>\epsilon$}{
\hspace{0.5cm}\textbf{solve} $\uu_k \in W$ satisfying
\begin{equation}
\begin{split}
a^\tau(\uu_k, \ww) =& F^n(\ww) - \int_{\Omega}(\sdiv{\vv_k})( \sdiv{\ww} ) dx \quad\forall \ww\in W, \\
\vv_{k+1} =& \vv_k + \rho \uu_k
\end{split}
\end{equation}}

\item[(2)] Set $\uu^{n+1} = \uu_k$ and define the pressure according to the Unified Stokes Algorithm \cite{lrsBIBia}:
\begin{equation}
  p^{n+1} = -\nu \Pi_P (\sdiv{\vv_k}), 
 \end{equation}
 where $\Pi_P$ is the $L^2(\Omega)$ projection onto $P= \nabla \cdot W$.
 \end{enumerate}
\end{algorithm}

To generate the mesh, we used {\tt mshr} \cite{lrsBIBih} with two inputs,
the meshsize $M$ and the number of segments $S$ used to approximate the circle.
The code was implemented in FEniCS \cite{alnaes2015fenics}.

For steady-state computations, we used Newton continuation \cite{lrsBIBih}
as a nonlinear solver for the steady-state Scott--Vogelius--Nitsche method
as explained in \cite{lrsBIBkj}.

\section{Drag and lift}\label{dragLift}

{The main focus of this paper is to analyze the dynamics of fluid flow around a cylinder in 2D.} Unfortunately, our solution space is infinite dimensional, so we may not directly observe the dynamics on this space. As such, we will instead use dimensionless coefficients associated to the drag and lift forces on the cylinder to define a two-dimensional phase space on which we may observe a projection of the dynamics of our system.

More precisely, {the \textit{drag force} is the component force on the cylinder} in the stream-wise direction. The \textit{lift force} is the component force on the cylinder acting in the direction perpendicular to the flow. For example, when air is flowing around the wing of an aircraft, the drag acts along the width of the wing, pointing towards the trailing edge. The lift acts opposite to the force of gravity, so it points vertically up.  In our simulations, we define the direction of lift to be along the positive vertical axis.

Drag and lift are typically identified by the dimensionless \textit{drag and lift coefficients},  $C_D$ and $C_L$, and they may each be further broken down into a \textit{pressure} component and a \textit{viscous} component \cite{lrsBIBjn}. The drag and lift coefficients corresponding to the pressure (a.k.a. pressure or form drag and lift) for a cylinder of radius 1 are given by
\begin{equation} \label{eqn:dragformpa} 
  C^d_P=- \oint_{\Gamma}  p \, \bn\cdot(1,0)\, ds, \qquad
  C^l_P=- \oint_{\Gamma}  p \, \bn\cdot(0,1)\, ds,
\end{equation}
where $\bn$ is the outward normal to $\Omega$, pointing into the cylinder on $\Gamma$, and $$\du=\half(\nabla u+\nabla u^t).$$

The viscous drag and lift (a.k.a. skin friction drag and lift) for a cylinder of radius 1 are given by 
\begin{equation} \label{eqn:dragformmb} 
C^d_V=\oint_{\Gamma} \nu \big((\du)(1,0)\big)\cdot\nn\, ds, \qquad
C^l_V=\oint_{\Gamma} \nu \big((\du)(0,1)\big)\cdot\nn\, ds.
\end{equation}
In our system, the pressure and viscosity are the only factors contributing to the drag and the lift. Therefore, the full drag and lift coefficients for a cylinder of radius 1 are given by 
\begin{equation} \label{eqn:dragfullmb} 
C^d=C^d_P+C^d_V,\qquad C^l=C^l_P+C^l_V.
\end{equation}

\begin{figure}[H] 
\centerline{\includegraphics[height=0.4\textwidth]{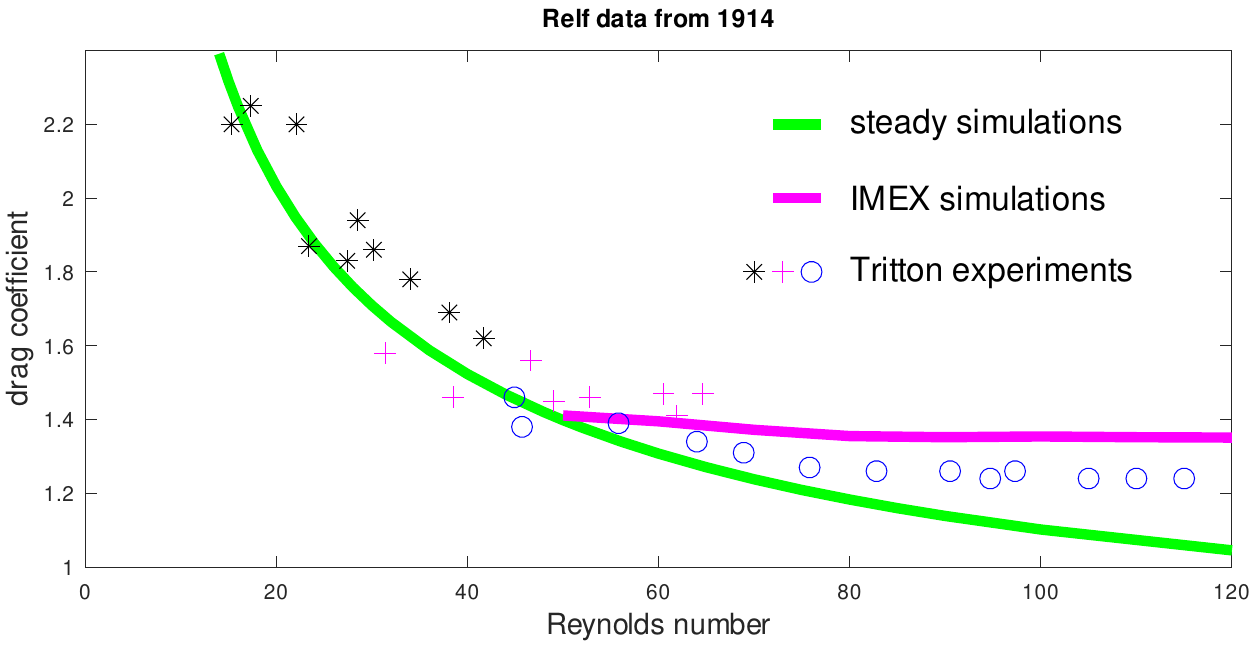}}
\caption{Data for cylinder drag from Tritton \cite{ref:trittonlowrexpts} together 
with simulations (solid lines) using the SVN scheme with
free-stream boundary conditions on the domain \eqref{eqn:domparam}
with parameters $L=300$, $b=30$, $w=30$. For the steady flow simulations, we set mesh resolution $M=128$, and segments =1024. For the time-dependent simulations, the time step was set to $\Delta t= 0.01$ with $M=32$ and Segments =2048.
}
\label{fig:trittonrelf} 
\end{figure}

Figure \ref{fig:trittonrelf} depicts computations of the \textit{mean} drag coefficient $C_D$ using our SVN scheme compared to the physical data collected in \cite{ref:trittonlowrexpts}.  Different symbols (and colors) indicate experiments with different 
cylinder size in \cite{ref:trittonlowrexpts}.
The longer curve represents simulations of steady flow, whereas the shorter (upper) curve represents 
time-dependent simulations.

\begin{figure}[H] 
\centerline{\includegraphics[height=0.45\textwidth]{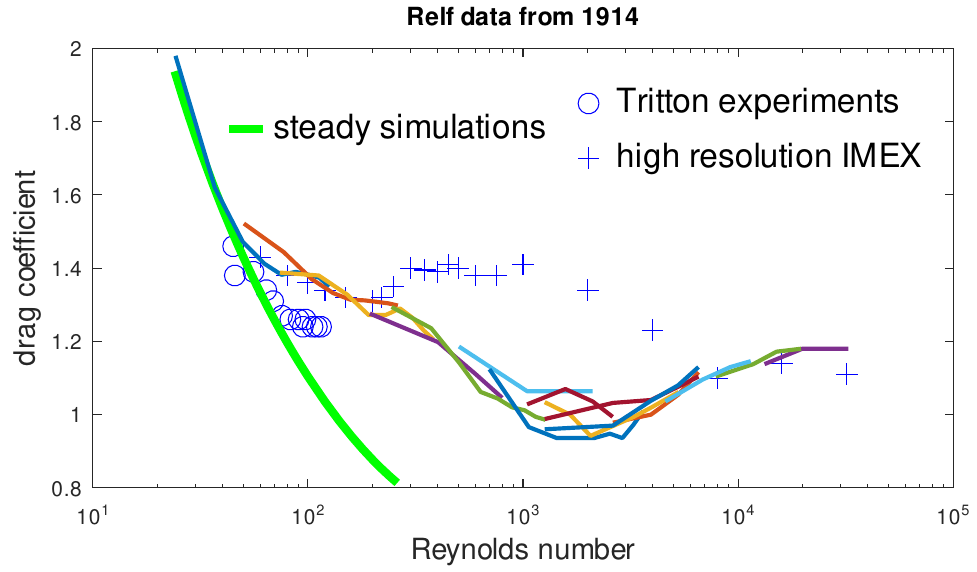}}
\caption{Data for cylinder drag compared to Relf \cite{ref:errelflowRex}. See Section \ref{sec:expData} for more details.
}
\label{fig:relfdata} 
\end{figure}

In Figure \ref{fig:relfdata}, we compare our numerical computations of the mean $C_D$ to the data from Relf \cite{ref:errelflowRex},  where the short, narrow, colored line segments indicated different wire diameters (see \ref{sec:expData}).  The wide red line indicates steady SVN simulation data and the magenta $\times$'s
indicate IMEX SVN simulations, both generated using the schemes
in Figure \ref{fig:trittonrelf} for channel dimensions given there.
The circles indicate the same Tritton data as shown in Figure \ref{fig:trittonrelf}.  Simulation parameters for the steady flow are the same as stated in Figure \ref{fig:trittonrelf}. For the time-dependent flow, we set $M=32$ and segments =2048,
$\Delta t=0.01$ for $ R< 300$,
$\Delta t=0.005$ for $300\leq R< 1000$,
$\Delta t=0.004$ for $R\geq 1000$.

\section{Dynamics perspective} \label{dyn}
For flow around a cylinder, it is well known that a Hopf bifurcation from steady, symmetric flow to the Karman vortex street occurs around $R = 50$ 
\cite{ref:cylinderHopfvortexstreet,chen1995bifurcation}. We observe this in our own data, as shown in Figures \ref{fig:trittonrelf} and \ref{fig:relfdata} where the data from the steady flow and time-dependent flow simulations diverge. At this point, a \textit{seemingly periodic} solution arises for the linearized, time-dependent flow.  As the Reynolds number continues to rise, however, this periodic solution will break down, to be replaced with more chaotic flow.

{The breakdown of this periodic flow and transition to chaotic flow is a process of great scientific interest, however it becomes challenging to study due to the complex dynamics of flow at these higher Reynolds numbers. It is therefore our goal to study these dynamics from a different perspective--\textit{numerical computation of the drag and lift coefficients, $C_D$ and $C_L$.} }We see from \eqref{eqn:dragformpa} and \eqref{eqn:dragformmb} that $C_D$ and $C_L$ are dependent on the solutions of the velocity and pressure, $(\uu,p)$, so if our flow is periodic, this periodicity should be reflected in the drag and lift coefficients as well.  

\begin{figure}[H] 
\vspace{-0cm}
\centerline{\includegraphics[height=1.6in]{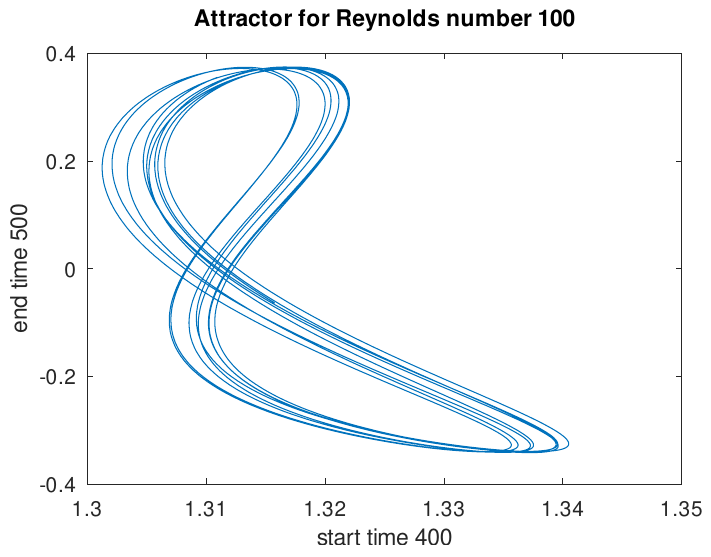}
            \includegraphics[height=1.6in]{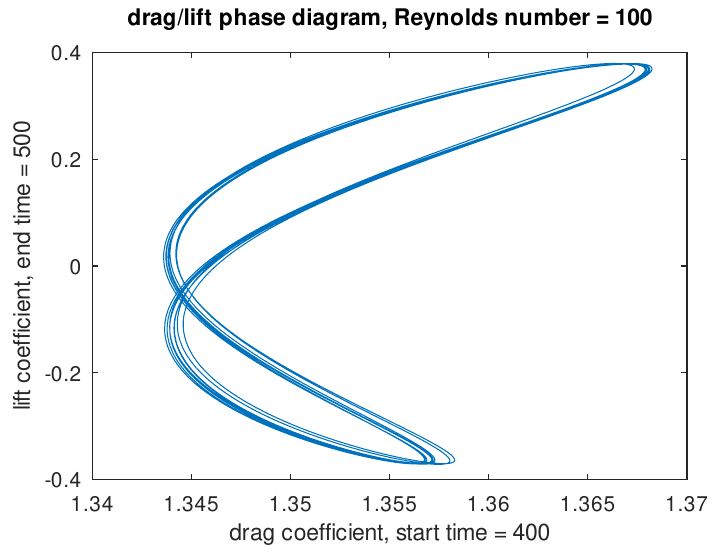}
            \includegraphics[height=1.6in]{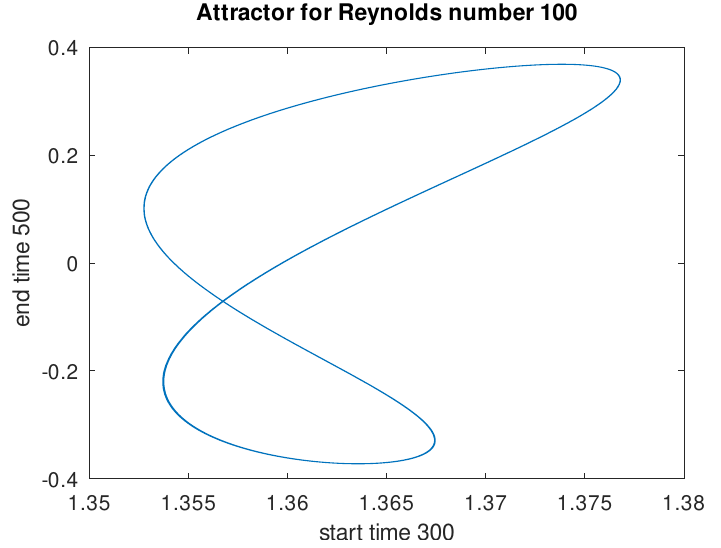}}
\vspace{-2mm}
\caption{Effect of computational resolution on attractor, $R=100$.
$M=8, 16, 32$ (left to right), with segments = $64M$.
$\Delta t=0.02$ except for $M=32$ where $\Delta t=0.01$.
Mean drag values were 1.32, 1.35, 1.36, respectively,
for the intervals indicated in each plot.}
\label{fig:vazmeshref} 
\end{figure}

In Figure \ref{fig:vazmeshref}, we plot the lift-versus-drag phase diagrams
for $R=100$ for increasingly refined discretization parameters.
More precisely, if the drag and lift values are $x_i$ and $y_i$ at times $t_i$
(we chose the values $t_i$ to be all of the time steps computed), then the
plot is the set of points $(x_i,y_i)$ connected by straight lines, with the
values of $i$ chosen so that $t_i$ is not smaller than the indicated start time.

\begin{figure}[H] 
\vspace{-0cm}
\centerline{\includegraphics[height=1.6in]{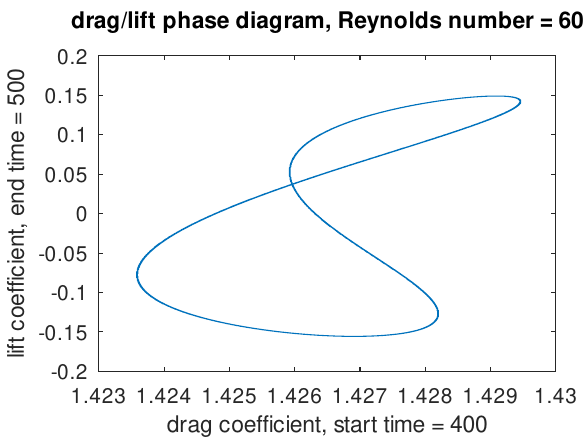}
            \includegraphics[height=1.6in]{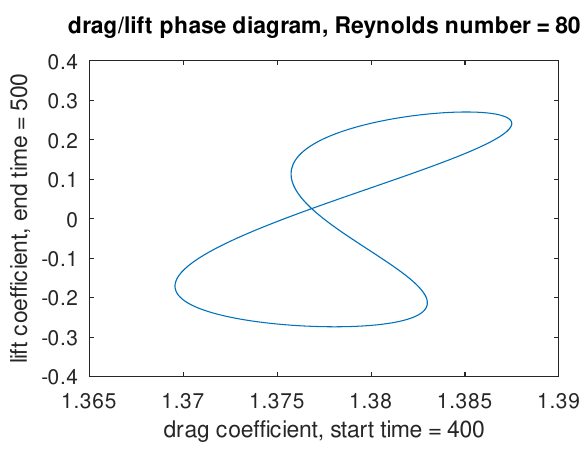}
            \includegraphics[height=1.6in]{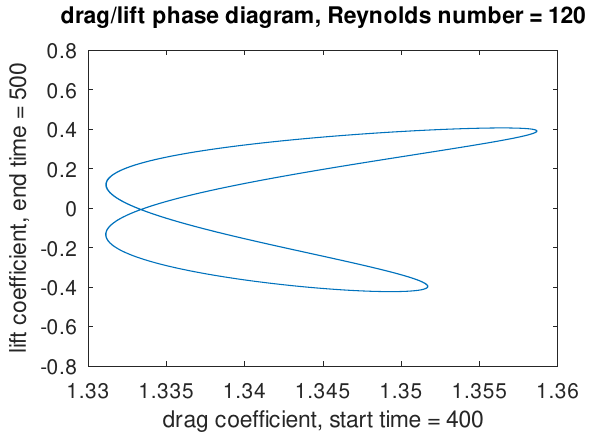}}
\vspace{-2mm}
\caption{Attractor evolution, $R=60, 80,120$. $M=32$,  segments = 2048, and $\Delta t=0.01$.
Mean drag values were 1.43, 1.38,  and 1.34, respectively,
for the intervals indicated in each plot.
See Figure \ref{fig:vazmeshref} for $R=100$.}
\label{fig:atraktr} 
\end{figure}

At this Reynold's number, the flow should still be largely periodic.
Indeed, what we see in the phase diagrams for $R=100$ with different numerical
resolution is a convergence of the phase diagram to one representing
periodic flow.
The geometry of the phase diagrams for the less-refined simulations are
similar, but it is clear that sufficient resolution is needed to draw
signficant conclusions. 
And yet, the drag and lift coefficients are undeniably periodic at $R=100$.
This periodic behavior is also seen at lower Reynold's numbers,  as shown
in Figure \ref{fig:atraktr}.

Figure \ref{fig:vazmeshref} shows an interesting feature of the periodic
vortex street, namely, that it is not up-down symmetric.
For $M=8$, drag is higher for negative lift, whereas for the larger messhes
it switches the other way around.  Such asymmetry is not entirely unexpected, particularly if we compare the results with a close qualitative assessment of well-known images capturing the vortex street, such as those in assembled in Milton Van Dyke's \textit{Album of Fluid Motion}\cite[Figures 95-97]{ref:vanDyke}. However, we cannot quantitatively verify this claim.

It is important to note that periodic solutions for $C_D$ and $C_L$ do not necessarily imply periodic $(\uu, p)$.  As such, we cannot directly make claims about the solution to \eqref{eqn:nastoensn} from these phase plots. However, we do expect there to be significant correlation with the solution. 

Furthermore, the dynamics of drag and lift are of significant interest on their own, as they may allow us to paint a more detailed picture of the flow around the cylinder and identify interesting features of that we may otherwise not observe. For example, we can observe in Figure \ref{fig:relfdata} that the mean drag for the oscillating flow is significantly
different from the drag for the steady, symmetric flow for $R>>50$. {In other words, at the onset of the Karman vortex street, the seemingly periodic solution that arises does not oscillate around the steady-state solution.  If this were the case, we would expect the time-dependent drag to oscillate around the steady-flow drag, resulting in the same mean. The fact that this does not occur suggests that the flow fundamentally changes at this point.}

\begin{remark}
We note that the behavior observed in these phase plots is subject to a \textit{spin-up time}.  More specifically, we have performed all of our studies by solving \eqref{eqn:nastoensn} starting
with $u(t=0)$ given by the Stokes solution for the boundary value problem.
Consequently, in the initial time development, the steady solution for the Navier--Stokes
system \eqref{eqn:nastoensn} for the given Reynolds number emerges,
and then (for Reynolds numbers 50 and above)
the vortex street begins to emerge, due to the asymmetry of the mesh.
(If the mesh were perfectly up-down symmetric, the steady Navier--Stokes
solution would persist.)
After further development (\textit{spin-up time}), the lift and drag increase, eventually reaching
a state of persistent oscillations.  We see this demonstrated in Figure \ref{fig:timevoluti}, as the phase plot at $R=200$ is significantly more periodic on the time-step interval $[400,500]$, whereas we still observe the solution pre-periodicity if we start at time-step $200$.

From our investigations, the time it takes for a solution to settle in to oscillations depends inversely on the Reynolds number. Therefore, in the investigations that follow, we provide information on the starting and ending time-steps whenever it is available, and relevant.

\begin{figure}[H] 
\vspace{-0cm}
\centerline{\includegraphics[height=2.0in]{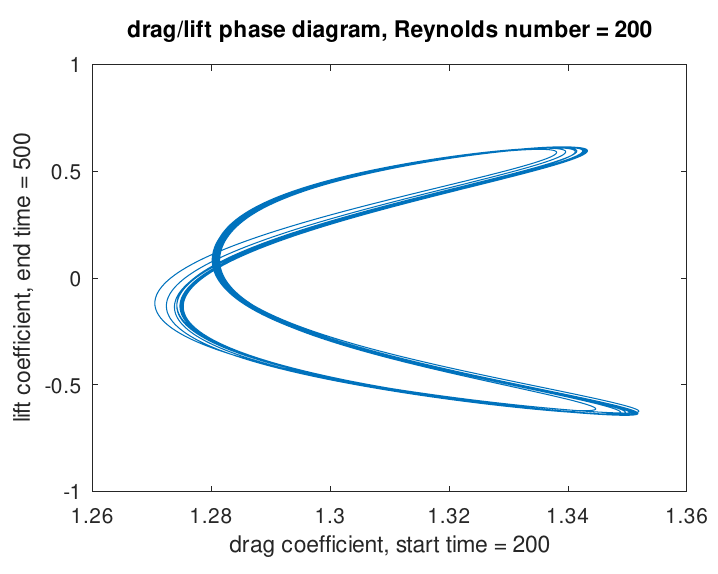}
            \includegraphics[height=2.0in]{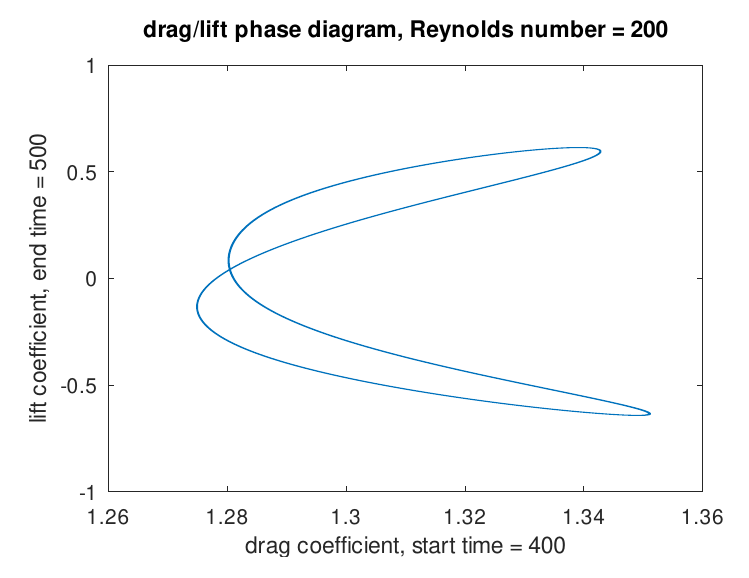}}
\vspace{-2mm}
\caption{Attractor evolution in time, $R=200$.
Left: start time is $T=200$.
Right: start time is $T=400$.
$M=32$, $\Delta t=0.01$, segments = 2048.
Mean drag values were 1.311 for both cases.
}
\label{fig:timevoluti} 
\end{figure}  
\end{remark}

{In this paper, we choose to specifically investigate the oscillatory behavior of the drag and lift as Reynolds number increases using a variety of mathematical tools, in order to more fully understand the development of chaotic flow. }

\begin{figure}[H] 
\vspace{-0cm}
\centerline{\includegraphics[height=1.5in]{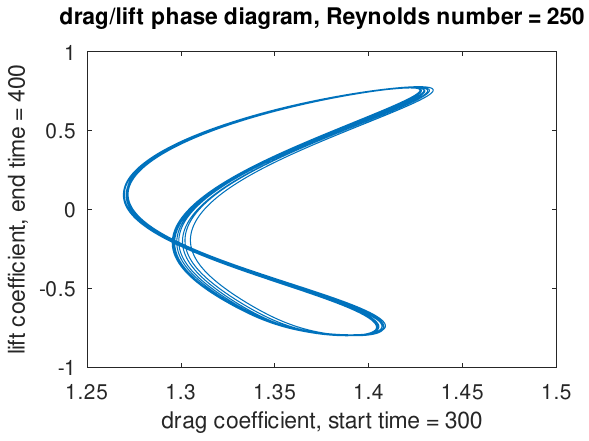}
            \includegraphics[height=1.5in]{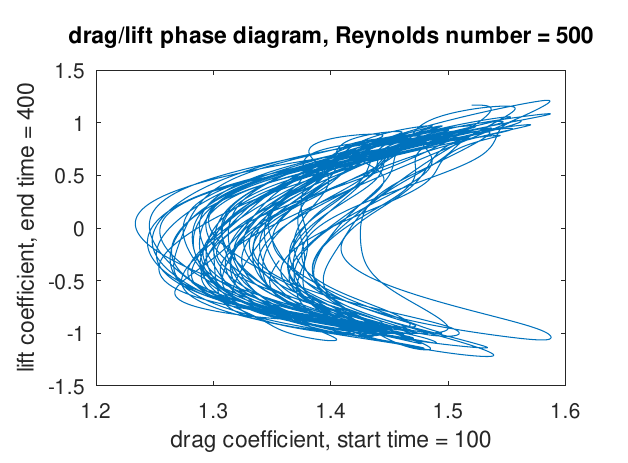}
            \includegraphics[height=1.5in]{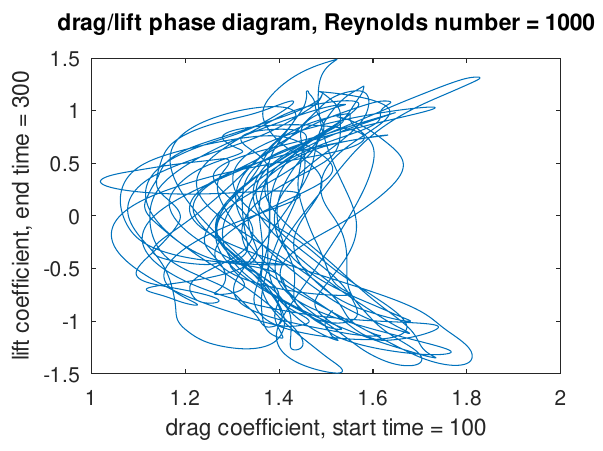}}
\vspace{-2mm}
\caption{Attractor evolution, $R=250, 500, 1000$. $M=32$, segments = 2048,
$\Delta t=0.005, 0.005, 0.004$.
Mean drag values were  1.3531 ($R=250$), 1.4027 ($R=500$), and 1.4143 ($R=1000$) 
for the time intervals indicated.
}
\label{fig:Furstchaos} 
\end{figure}

A preliminary investigation of the phase plots for $R=250, 500,1000$ as shown in Figure \ref{fig:Furstchaos} demonstrates the breakdown of periodicity that we expect at higher Reynolds numbers.  However, the plots at $R=250$ and $R=500$ indicate that this chaotic behavior evolves over time. {As such, the main goal of the proceeding sections will be to \textit{quantify} the behavior observed in the plots, so that we may track the evolution and onset of chaotic flow.}

\begin{remark}
We note that the clear chaos at $R=1000$ appears to contradict existing literature.  In particular, we note that \cite[Figure 5]{ref:wolfchaocylinderflow} suggests that drag
is nearly periodic at $R=1000$, not chaotic.

The discrepancy may be due to computational dissipation in the simulations
in \cite{ref:wolfchaocylinderflow}.
To make a direct comparison, we plot our values for the lift as a function 
of time in Figure \ref{fig:plotlift}, which should be compared 
to \cite[Figure 3]{ref:wolfchaocylinderflow}.
The scheme used in \cite{ref:wolfchaocylinderflow} is the
vortex blob method of Chorin.
The numerical diffusion and dispersion associated with this method 
has been examined in \cite{oliver2001vortex}.
\end{remark}

\subsection{Lyapunov exponent}

Lyapunov exponents \cite{wiggins2003introduction} are used to quantify the chaotic behavior of a system. As explained in \cite{ref:timeseriesLyapunov}, if a system contains at least one positive Lyapunov exponent, it is considered to be chaotic. The magnitude of the exponent gives an indication of how quickly this chaotic behavior appears. Thus, a large, positive value indicates a highly chaotic system with very little spin-up time.

A description of the method used to determine the Lyapunov exponents is given in Appendix \ref{sec:mleo1}.  To summarize,  we {begin the method by computing a \textit{best approximate period} for our data, $p_S$ (see Appendix \ref{avP}). Then,  we consider a vector of data $x(t_0)$, called a \textit{delay coordinate}, covering $M$ time steps,  $\{t_0, t_1, \cdots, t_M\}$, separated by a time parameter $\Delta$, so that $t_k = t_0 + k \Delta$. }

{Then, beginning at time $t_0$, we take a segment $x(t_0)$ of our drag and lift data consisting of the values at $m$ sequential time steps and find a later segment of data, $x(t_0')$ that is the \textit{nearest neighbor} (closest in value) to our initial segment occurring at a time $t_0'$ close to $t_0 + p_S$.  We repeat this process for each $t_k$, and for each iteration we compute the logarithm (base 2) of the ratio of the distances between nearest neighbors at their original times $t_k$ and $t_k'$ and at the evolved times $t_k + \Delta$ and $t_k' + \Delta$. The first Lyapunov exponent, $\lambda_1$, is then approximated by taking the average of these quantities over the time intervale $[t_0, t_M]$.}

\begin{figure}[H] 
\vspace{-0cm}
\centering
\includegraphics[height=2.2in]{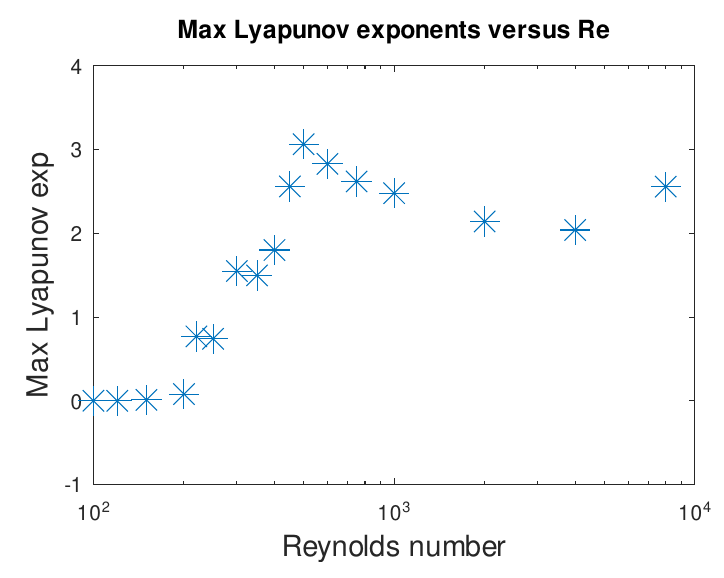}
\vspace{-1mm}
\caption{ Maximal Lyapunov exponent $\lambda_1$ as computed by the algorithm.
For $R\leq 200$, data from the time interval $[380,480]$ was used; 
for $220\leq R\leq 600$, data from the time interval $[280,380]$ was used;
and for Reynolds number 1000, data from the time interval $[180,280]$ was used.
For the simulations, meshsize $M=32$, segments=2048, time steps 
$\Delta  t=0.1$ for $R\leq 250$,
$\Delta  t=0.05$ for $300\leq R\leq 600$, and
$\Delta  t=0.04$ for $R=1000$.
}
\label{fig:mleo} 
\end{figure}

Figure \ref{fig:mleo} indicates how
the maximal Lyapunov exponent $\lambda_1$ grows with the Reynolds number $R$.
The data suggests that the flow becomes increasingly chaotic
as $R$ increases, although a maximum is reached at $R=500$.  This may be a result of numerical error, or it may indicate that some sort of bound for chaotic behavior. 

{Of particular note is how these results confirm an informal visual conclusion regarding the onset of chaotic behavior. We note from Figure \ref{fig:timevoluti} that the phase plots at $R=200$ suggest that the flow is periodic, yet by $R=250$ we observe the onset of chaotic flow, as indicated by Figure \ref{fig:Furstchaos}. This is confirmed in Figure \ref{fig:mleo} with the first noticeably nonzero value of $\lambda_1$ appearing just after $R=200$.}


\subsection{Strouhal correlation}
Lyapunov exponents provide a guide for us to quantify the onset of chaotic flow, and they appear to indicate through our drag and lift computations that our system develops chaotic behaviour gradually after the onset of the Karm\'an vortex street. As the onset of trubulence indicates the breakdown of oscillatory flow, a question arises as a natural companion to our Lyapunov studes: how \textit{periodic} is the flow as Reynolds number increases?

\begin{figure}[H] 
\vspace{-0cm}
\centerline{(a)\includegraphics[height=2.2in]{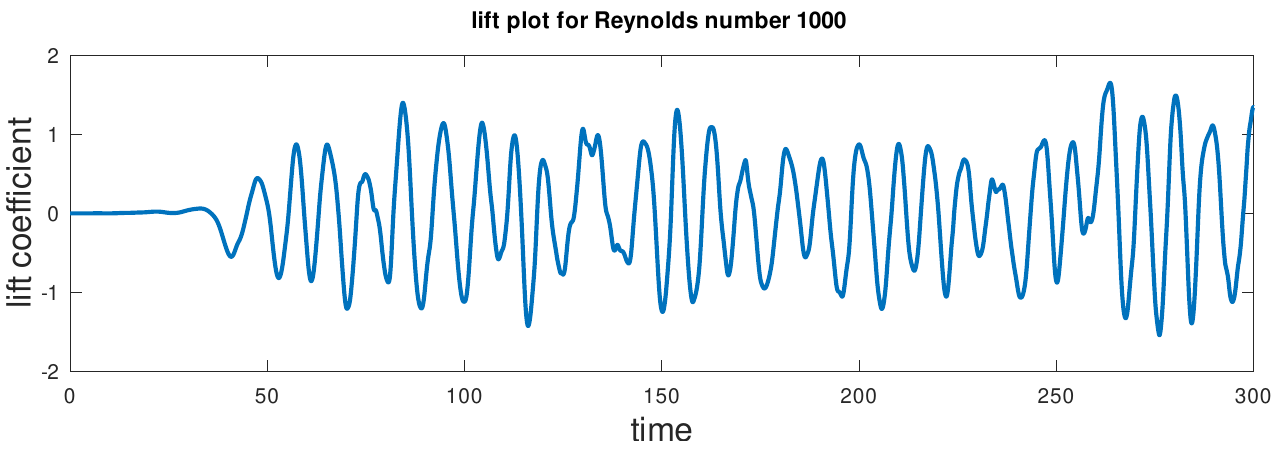}}
\centerline{(b)\includegraphics[height=2.2in]{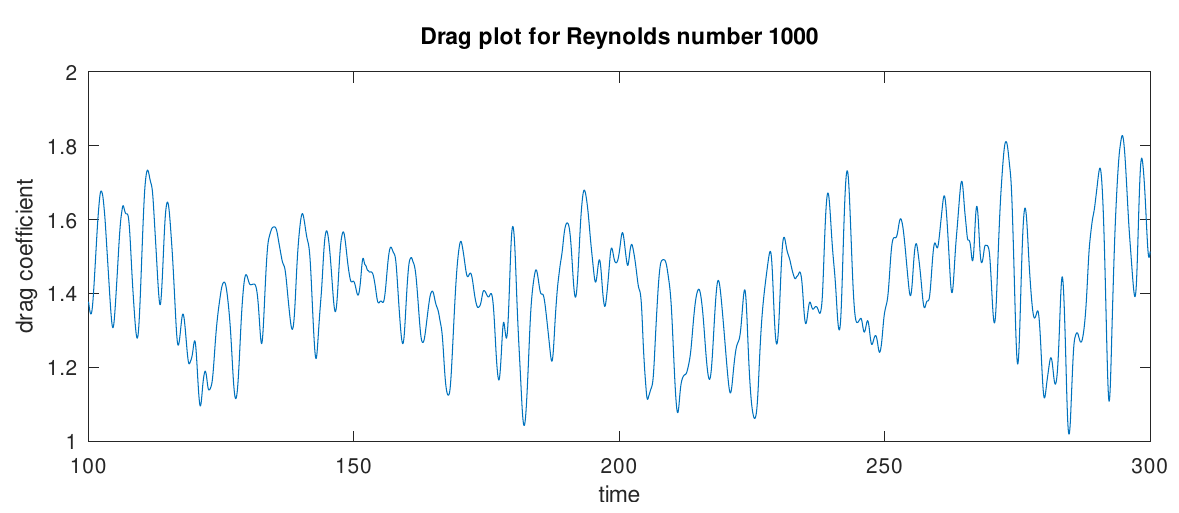}}
\vspace{-2mm}
\caption{Lift (a) and drag (b) as functions of time for $R=1000$ for $t\in[0,300]$.
See Figure \ref{fig:mleo} for more details on the simulation parameters.}
\label{fig:plotlift} 
\end{figure}

One method to quantify the periodicity of the drag or lift is to simply view the data as a function of time, as depicted for $R=1000$ in Figure \ref{fig:plotlift}, and measure how close the function is to a known periodic function.  To do this, we may compute the \textit{Strouhal 
period} $\pi_S$ \cite{birkhoff1953formation} for a given Reynolds number $R$
by fitting $\sin(kt)$ to the drag or lift data.
Typically, the lift is used.

We vary $k$ until the correlation $c_S$ between the sinusoidal function and the
lift data is a maximum.
More precisely, we define 
\begin{equation} \label{eqn:stroucorr} 
c_S=\frac{L^t \Sigma}{L^t L},
\end{equation}
where $L$ is the lift data and $\Sigma$ is the sinusoidal data on the same time grid. A value of $c_s$ close to $1$ indicates a high correlation between the lift data and a periodic function. In other words, if $c_s \approx 1$, the data is nearly periodic.

 We see from the data in Table \ref{tabl:stroohall} that it is eventually
difficult to define the Strouhal data as the Reynolds number increases, as 
the correlation between the lift data and a periodic function becomes negligible. {However, for lower Reynolds number, we do see values of $c_s$ close to $1$, matching the periodic behavior we observe in the phase plots.} As Reynolds number increases, we clearly see that $c_s$ gets farther from $1$, indicating a breakdown of periodicity. This is gradual at first, as $c_s$ has only decreased to $0.972$ by $R = 200$. However, by $R=1000$ it has dropped all the way to $0.087$.

Additionally, we note that $\pi_S$ is not the same as the fixed point $p_S$ found when estimating the Lyapunov exponent, and the divergence of the two quantities is an interesting aspect to consider.  As we can see in Table \ref{tabl:periods}, the computed values of $\pi_S$ and $p_S$ are very similar for lower Reynolds numbers, corresponding to the range where we expect nearly periodic behavior, as indicated by values of $c_S$ close to $1$. However, as the Reynolds number increases, the value of $\pi_S$ begins to decline rapidly, corresponding to a similar decline in $c_S$ (see, for reference, the values for $R=450$ and $R=1000$). This appears to indicate a decline in \textit{periodic} behavior of the solution as $R$ increases. The value of $p_S$, however, remains more stable, although it does decline as $R$ increases. {One possible interpretation of this is that the solutions retain some oscillatory behavior as Reynolds number increases, however they are not periodic due to the onset of chaotic flow.}

\begin{table}
\begin{center} 
\begin{tabular}{|c||c|c|c|c|c|c|}
\hline
$R$ & $\pi_S$  & $p_S$   & $c_S$  & start time & end time & $\Delta t$ \\
\hline
    60 & 14.43  & 14.58 & 0.967 & 400 & 500& 0.01 \\
    80 & 13.13  & 13.14  & 0.992 & 300 & 500& 0.01 \\
   100 & 12.09  & 12.15 & 0.983 & 300 & 500& 0.01 \\
   120 &  11.53  & 11.52  & 0.980 & 300 & 500& 0.01 \\
   150 & 11.02  & 10.96  & 0.975 & 300 & 500& 0.01 \\
   200 &  10.66  & 10.64  & 0.972 & 300 & 500& 0.01 \\
   450 &  8.90  &  9.05 & 0.908 & 200 & 400& 0.005 \\
  1000 &  6.99  &  8.79  & 0.087 & 100 & 300& 0.004 \\
\hline
\end{tabular}
\end{center}
\vspace{-6mm}
\caption{Strouhal data. See Figure \ref{fig:mleo} for simulation details.  To compute $p_S$, we used $m=5$, and the individual $p$ values were computed
at intervals of $t=0.125$.}
\label{tabl:periods}
\end{table}

\subsection{Strouhal Number}\label{strouhal}
There are other numbers of interest in addition to the periods $p_S$ and $\pi_S$ and
the Strouhal correlation $c_S$. For instance,  the \textit{shedding frequency}, $f_S$, may be defined by $f_S=1/p_S$ or $f_S=1/\pi_S$.  This quantity should reflect the frequency of vortex shedding from the trailing edge of the cylinder at the onset of the Karm\'an vortex street. As $\pi_S$ is determined to measure \textit{periodicity} of the flow dynamics, it is unlikely to be a useful measure for the rate of vortex shedding as $R$ increases and flow becomes more chaotic. As such,  we will define $f_S$ using the quantity $p_S$, as it is a more reliable measure of the shedding period at higher Reynolds numbers. The values of $p_S$ for higher $R$ are shown in Table \ref{tabl:stroomore}. These were computed in the same way as stated in Table \ref{tabl:periods}.

\begin{table}
\begin{center} 

\begin{tabular}{|c|c|c|c|c|}
\hline
 $R$ &  $p_S$   &  start time & end time & $\Delta t$ \\
\hline
2000 &     10.77  &       180   &      380   &   0.002 \\
4000 &     12.07  &       280   &      480   &   0.002 \\
8000 &    11.02  &       180   &      380   &   0.002 \\
\hline
\end{tabular}
\end{center}
\vspace{-6mm}
\caption{
See Figure \ref{fig:mleo} for simulation details.  Shedding frequencies at higher Reynolds numbers.}
\label{tabl:stroomore}
\end{table}

From the shedding frequency, the \textit{Strouhal number} $N_S$ is defined to be $f_S\ell/U$, where 
$\ell$ and $U$ are the characteristic length and speed used
to define the Reynolds number, so $\ell=2$ and $U=1$ in our case. Thus $N_S$ is understood to be the shedding frequency of the flow around the cylinder, scaled by the characteristic length and speed.  Additionally, the Strouhal number is closely related to the inverse of the
Keulegan--Carpenter number \cite{keulegan1956forces}.

In \cite{ref:FeyStrouhalFormula}, an empirical estimate, $\phi_S$, for the Strouhal number was provided as a piecewise-linear function of the Reynolds number for $47 < R < 2\times 10^5$,
\begin{equation}\label{empirical}
\phi_S(R) = S_r^* + \frac{m}{\sqrt{R}},
\end{equation}
where the quantities $S_r^*$ and $m$ depend on the Reynolds number (see \cite[Table 1]{ref:FeyStrouhalFormula}). These results indicate that the shedding frequency should increase until $R>1300$, at which point it should begin to decrease. 

\begin{table}
\begin{center} 

\begin{tabular}{|c|c|c||c|c|c|}
\hline
$R$ &  $N_S$ & $\phi_S$  & start time & end time & $\Delta t$ \\
\hline
    60 & 0.1372 & 0.1347  & 400 & 500& 0.01 \\
    80 & 0.1522 & 0.1526 & 300 & 500& 0.01 \\
   100 & 0.1646 & 0.1648  & 300 & 500& 0.01 \\
   120 & 0.1736 & 0.1739 & 300 & 500& 0.01 \\
   150 & 0.1815 & 0.1838  & 300 & 500& 0.01 \\
   200 & 0.1880 & 0.1828  & 300 & 500& 0.01 \\
   450 & 0.2237 & 0.2049  & 200 & 400& 0.005 \\
  1000 & 0.2276 & 0.2118 & 100 & 300& 0.004 \\
  2000 & 0.1857 &   0.2116  &       180   &      380   &   0.002 \\
4000 & 0.1657 &    0.2093  &       280   &      480   &   0.002 \\
8000 & 0.1815 &    0.2022  &       180   &      380   &   0.002 \\
\hline
\end{tabular}
\end{center}
\vspace{-6mm}
\caption{
Strouhal data.
Simulation data is the same as Figure \ref{fig:relfdata}.}
\label{tabl:stroohall}
\end{table}

 As we can see from Table \ref{tabl:stroohall}, our computed values of $N_S$ appear to agree nicely with the empirical quantity $\phi_S$ for lower Reynolds numbers. The agreement is less significant at higher Reynolds numbers, however we note that this is likely due to computational error in the approximation of $p_S$. It is, however, significant to note that the behavior of $N_S$ with respect to $R$ is indeed what we expect from the empirical results in \cite{ref:FeyStrouhalFormula}. Indeed, we observe $N_S$ increasing with $R$ until we surpass $R=1000$. At this point, the shedding frequency decreases with $R$, much like $\phi_S$.  {It is noted in \cite{ref:FeyStrouhalFormula} that this change corresponds to the known onset of a Kelvin Helmholtz instability in the separated shear layer at $R \approx 1300$. }

\subsection{Fractal dimension}
The final method used to investigate the flow dynamics is the fractal dimension. The fractal dimension $D$ is defined by \cite[(3)]{ref:fracordimeapp} 
in terms of the proximity of pairs of drag/lift data.
Let $x_1(t_i)$ and $x_2(t_i)$ denote the drag and lift, respectively,
at a time step $t_i$, and let $\xx(t_i) = (x_1(t_i),x_2(t_i))$.

Define $C(r)$ by
$$
C(r)=\#\set{(i,j)}{j>i,\;|\xx(t_i)-\xx(t_j)|\leq r,\;\tau\leq t_i,t_j\leq T},
$$
where $\#$ means the cardinality of the set (the number of pairs),
$\tau$ is a start time chosen to omit the start-up phase of the
simulations, and $T$ is the total length of the simulation.
Then $D$ is defined by fitting the expression
$$
C(r)\approx c r^D,
$$
for $c$ a fixed constant, then solving for $D$.

One way to do this is by computing $C(r)$ for 
$r=r_0,2r_0,4r_0,\dots,2^kr_0$ for some integer $k$.

We then expect that
$$
\frac{C(r)}{C(r/2)}=2^D.
$$
Thus for each $r$ we can define $D_r$ by
$$
D(r)=\log\big(C(r)/C(r/2)\big)/\log 2.
$$

Essentially, $D(r)$ is a measure of the global spread of data. If data is dense, then there will be a large number of data points clustered within a distance $r$ of any given point. Thus, we expect the value of $D$ to decrease as data becomes more chaotic and, consequently, more spread out.

\begin{figure}[H]
\centering
\includegraphics[height=2.2in]{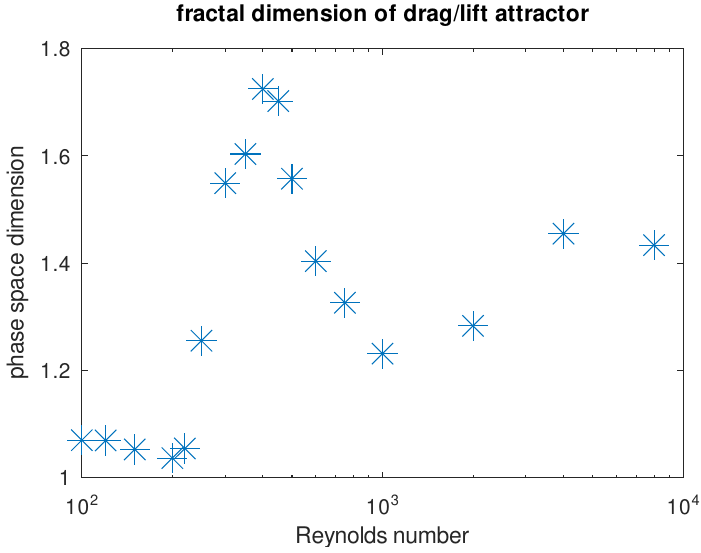}
\caption{Fractal dimension of attractors, the transition zone.
For $R\leq 200$, the time interval $[400,500]$ was used, 
for $220\leq R\leq 600$, the time interval $[300,400]$ was used, 
and for Reynolds number 1000, the time interval $[200,300]$ was used.
The sampling intervals were $\Delta  t$ for all cases.  For the simulations, meshsize $M=32$, segments=2048, time steps 
$\Delta  t=0.1$ for $R\leq 250$,
$\Delta  t=0.05$ for $300\leq R\leq 600$, and
$\Delta  t=0.04$ for $R=1000$.
}\label{fractalDim}
\end{figure}

Table \ref{tabl:fradimex} shows typical data for one value of
Reynolds number ($R=500$) where  $r_0=0.00125$ and $k=6$. 
We see that the dimension $D(r)$ varies, but is fairly stable in
the middle range of $r$ values.
To be consistent, we took $D(0.01)$ as the representative value
for all values of $R$.
Figure \ref{fractalDim} shows the evolution of the fractal dimension $D(0.01)$ 
of the attractors for a range of Reynolds numbers.


\begin{table}
\begin{center} 
\begin{tabular}{|c|c|c|}
\hline
$r$ & $C(r)$  & $D(r)$   \\
\hline
     0.00125&       19013 &    NA       \\
     0.0025 &      53436  &    1.4908    \\
      0.005 & 1.5419e+05  &    1.5289     \\
       0.01 & 4.5383e+05  &    1.5574      \\
       0.02 & 1.4391e+06  &     1.665       \\
       0.04 & 4.6319e+06  &    1.6864        \\
       0.08 & 1.3105e+07  &    1.5005         \\
\hline
\end{tabular}
\end{center}
\vspace{-6mm}
\caption{Fractal dimension data for $R=500$.}
\label{tabl:fradimex}
\end{table}

\color{black}

As shown in Figure \ref{fractalDim}, we observe low values for $D$ when $R$ is low and the phase plots indicate more periodic behavior. As $R$ increases, however, so does the value of $D$, indicating that our flow is no longer periodic. 

A particularly interesting feature of the fractal dimension data is that the value of $D$ eventually peaks and begins to decrease before reaching $R=1000$. The increase of $D$ and subesequent decrease correlates closely with
the same behavior of the Lyapunov exponent.

\section{Data availability statement}

The simulation  data, the codes used to generate it, and the analysis codes
will all be posted on a sutiable website, such as Zenodo, once the paper is accepted
for publication.

\section{Conclusions}

In this paper, we provide a computational study of the flow dynamics for Reynolds numbers up to $\approx10^4$. Our results provide quantitative evidence that the vortex shedding in the Karm\'an vortext street is periodic for Reynolds number in the approximate range of $50$ to $250$, at which point the flow becomes aperiodic. Additionally, we provide evidence that the time average of oscillatory flow that arises at the onset of the Karm\'an vortex street is notably different than the steady flow. 

We do, however, note that there is a discrepancy between our computational results and experimental data. This may be due to phenomena in the physical experiments that we have not yet investigated, such as vibrations in the cylinder. Future directions for this work therefore include a computational investigation of the impact of these vibrations on the flow dynamics. We are also interested in studying flow dynamics using different computational methods.

\appendix
\section{Lyapunov exponent}
\label{sec:mleo1}

\begin{figure}[H] 
\vspace{-0cm}
\centerline{\includegraphics[height=4.0in]{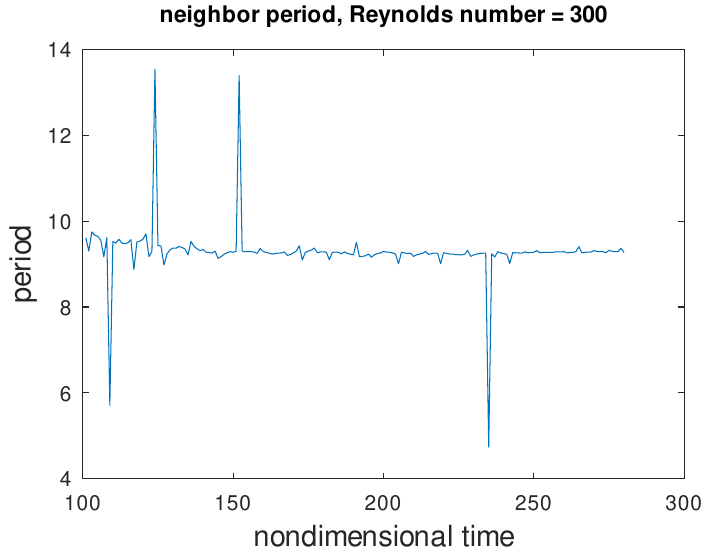}}
\vspace{-2mm}
\caption{Period evaluation for Re=300.  Simulation data is the same as Figure \ref{fig:relfdata}.}
\label{fig:example300} 
\end{figure}

Lyapunov exponents \cite{wiggins2003introduction} are key 
metrics for describing dynamic behavior.
They are inherently local metrics, evaluating the local behavior of
trajectories.
Their mathematical definition involves a limit as time goes to infinity, 
but a popular method \cite{ref:timeseriesLyapunov} provides a way 
to approximate the largest exponents with finite data, by
measuring aggregate local behavior of the dynamical system.
First the local spread of nearby orbits is measured throughout the system, and
then this data is aggregated over the entire dynamical data set. {For this paper, we adapted the method described in \cite{ref:timeseriesLyapunov} to suit our simulated data. In particular, the authors in \cite{ref:timeseriesLyapunov}, do not find an approximate period to restrict their search when identifying nearest neighbors.  We however, choose to assume that our data is ``nearly-periodic'' so that we can limit the region in which we search for nearest neighbors, making the process more efficient. The result, which we present in this appendix, are two methods adapted from \cite{ref:timeseriesLyapunov}: one for finding an approximate \textit{Strouhal period} (which is necessary for the analysis in Section \ref{strouhal}), and one for approximating the maximal Lyapunov exponent, $\lambda_1$.} 

There are three parameters for the method \cite{ref:timeseriesLyapunov}.
An integer $m\geq 1$ is used to create ``delay cooordinates''
which creates higher-dimensional data.
A second (positive) paramter $\tau$ defines the delay distance.
It can be ``almost arbitrarily chosen'' \cite[section 5.1]{ref:timeseriesLyapunov},
so we choose it to be the time step $\Delta t$ for the simulations,
for simplicity.
Thus the delay coordinate data here is
$$
x_i=x(i\Delta t)
=[d(i\Delta t), \ell(i\Delta t),\dots,d((i+m-1)\Delta t), \ell((i+m-1)\Delta t)],
$$
where $d$ and $\ell$ are the drag and lift data.
Thus the delay coordinates correspond to short segments of the drag/lift
dynamics.

In describing the algorithm, we will pretend that $x$ is a continuous function
of time $t$ with values in $\Reyuls^{2m}$, but in the computations it is
treated as discrete time values.

The gist of the method in \cite[section 5.1]{ref:timeseriesLyapunov} is to
find two nearby (``nearest neighbor'') $x$ points, corresponding to
two different times, and to
consider the evolution of the dynamics from these two points, measuring 
the change in the distance $L$ between the evolved points.
Thus we seek nearby branches in phase space and measure the evolution
(increase or decrease) of distance between the two branches.

\subsection{Finding the period}\label{avP}

The first step in our approach to approximating the Lyapunov exponent is to find an approximate period for the data. To describe the method of finding the period, we assume that we have nearly periodic data, with period $p$. Then,
given an initial time point $t_0$, we find the smallest time $t_0'>>t_0$ for which
$x(t_0)$ and $x(t_0')$ are closest, see \cite[Figure 4a]{ref:timeseriesLyapunov} or
\cite[Figure 6]{ref:SkokosreviewArxiv08}.
More precisely, if $p$ is an initial approximation of the period of the data, we want to find $x(t_0')$, the \textit{nearest neighbor} to $x(t_0)$ that occurs at a time close to $t_0+p$. Therefore we want $t_0'-t_0\approx p$,
and not $t_0'-t_0\approx 0$ or $t_0'-t_0\approx 2p$. However, as this $p$ is an initial approximation to the period, the resulting $t_0'$ will not occur exactly at $t_0 + p$. Thus we define a mapping $\Phi(t_0,p)$ such that
\begin{equation}\label{eqn:fimap}
t_0'=\Phi(t_0,p),
\end{equation}
and $t_0'$ is the time at which the nearest neighbor of $x(t_0)$ occurs within an interval around $t_0+p$.

As a result, {the quantitity $p_0=t_0'-t_0 $ gives a refined estimate of the period, which we store before repeating the process at a later time.} Note $t'-t\approx p$ puts us on a different branch of our phase space. If we choose a different $p$, such as $t'-t=\Delta t$, we would be on the same branch, and the method would not work.

{In practice, we define $t_0'$ by minimizing the Euclidean distance $|x(t_0)-x(t')|$ over the interval
$t'\in [t_0 +p/2, t_0+3p/2]$.  In other words,
\begin{equation} \label{argmin}
t_0' = \rm{argmin}\{\vert x(t_0) - x(t')\vert : t' \in  [t_0 +p/2, t_0+3p/2]\}.
\end{equation}
Once $t_0'$ has been found, we} {we store the resulting value $p_0 = t_0'-t_0 = \Phi(t_0,p) - t_0$. We then move to a new delay coordinate some time $\sigma$ away from $x(t_0)$, say $x(t_1) = x(t_0 + \sigma)$,  and repeat the process by finding $t_1'$ and refining the period again to determine $p_1 = t_1'- \Phi(t_1, p)$.  We refer to Table \ref{tabl:fixptem} for examples of the effect of $\sigma$ on the computed period $p$. Subsequently, this nearest neighbor process continues for some $K$ iterations, at which point we run out of data segments of length $m$. Note that the base value of $p$ does not change for each iteration, but a refinement, $p_i$, is stored each time.} 

{When this process is completed, we take the average of the stored $\{p_i\}_{i=1}^k$ to define a new approximate period $\tilde{p}$. Subsequently, we may repeat the entire iterative process again, using $\tilde{p}$ in place of $p$.  In our research, we repeated this nearest neighbor iteration-averaging process a number of times until we reached an iteration in which the average of the $p_i$ values converged to a fixed point, which we label $p_S$. We take this fixed point to be our best approximation to the period of the data.}

Note that in some cases {during the iteration process}, we obtain values which are quite close to $p/2$ or
$3p/2$, which essentially indicates failure of the algorithm \eqref{argmin}{, as it suggests the current approximation of $p$ is either too large or too small to match the data}.
This is indicated in Figure \ref{fig:example300}, where there are
two large spikes up and two down.
The other smaller spikes indicate small variations in the estimate of $p$
as one moves along the trajectory.

\begin{figure}[H] 
\vspace{-0cm}
\centerline{\includegraphics[height=3.4in]{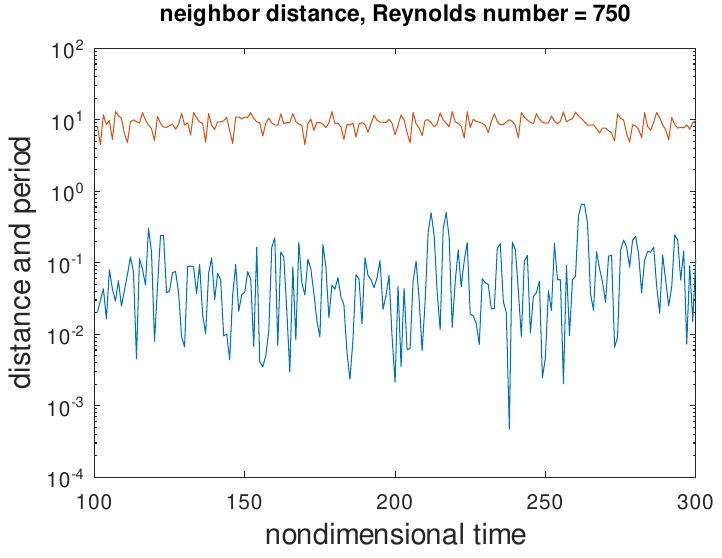}}
\vspace{-2mm}
\caption{Near neighbor distance $L$ for Reynolds number 750, in blue.
The red curve is the computed approximate period $p(t)=t'-t$,
with the fixed point (average) period being $p=8.9630$. Simulation data is the same as Figure \ref{fig:relfdata}.
}
\label{fig:nearnbraos} 
\end{figure}

This failure to locate a proper period does not happen for slower flows
(lower Reynolds numbers), but for faster flows it becomes more common, as shown in Figure \ref{fig:nearnbraos}.
{However, the subsequent averaging process appears to ameliorate these defects,} {and we note that we were able to determine the fixed point $p_S$ up to a prescribed accuracy for all Reynolds numbers.} A more sophisticated algorithm could improve this approach.

\begin{table}
\begin{center} 
\begin{tabular}{|c|c|c|c|c|c|}
\hline
  $R$  &    period  & t start   &   t end    &    $m$       & spacing $\sigma$ \\
\hline
          500   &   8.6036   &      280   &      380   &        5 &          1 \\
          500   &   8.776    &      280   &      380   &        5 &        0.5 \\
          500   &   8.8026   &      280   &      380   &        5 &       0.25 \\
          500   &   8.7928   &      280   &      380   &        5 &      0.125 \\
          500   &   8.7919   &      280   &      380   &        5 &     0.0625 \\
\hline
  750  &    8.8042  &       280 &        380   &        5   &        1\\
  750  &    8.9631  &       280 &        380   &        5   &      0.5\\
  750  &    8.9837  &       280 &        380   &        5   &     0.25\\
  750  &    8.9299  &       280 &        380   &        5   &    0.125\\
\hline
  750  &    8.9489  &       180 &        380   &        5   &     0.25\\
  750  &    8.9138  &       180 &        380   &        5   &    0.125\\
\hline
         1000   &   8.9035    &      80   &      280    &       5   &        1       \\
         1000   &   8.8135    &      80   &      280    &       5   &      0.5       \\
         1000   &   8.7546    &      80   &      280    &       5   &     0.25       \\
         1000   &   8.786     &      80   &      280    &       5   &    0.125 \\
         1000   &   8.7872    &      80   &      280    &       5   &   0.0625 \\
\hline
\end{tabular}
\end{center}
\vspace{-6mm}
\caption{Effect of the spacing parameter $\sigma$ on the computed period.
Simulation data is the same as Figure \ref{fig:relfdata}.}
\label{tabl:fixptem}
\end{table}

\subsection{Distances between orbits}

\begin{figure}[H] 
\vspace{-0cm}
\centerline{\includegraphics[height=2.4in]{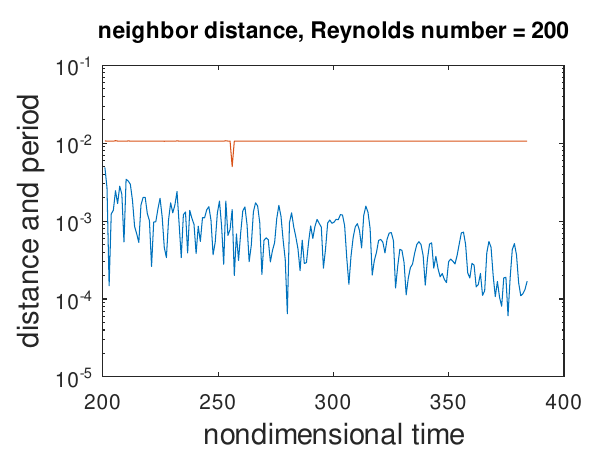}
            \includegraphics[height=2.4in]{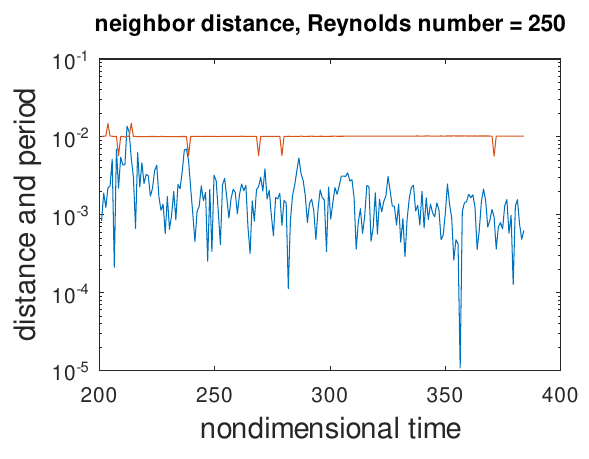}}
\centerline{\includegraphics[height=2.4in]{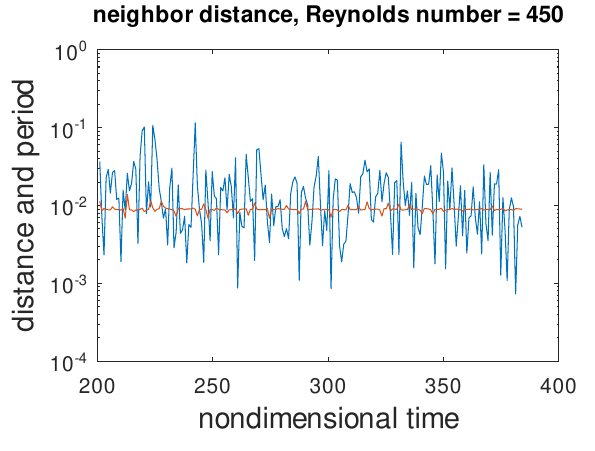}
            \includegraphics[height=2.4in]{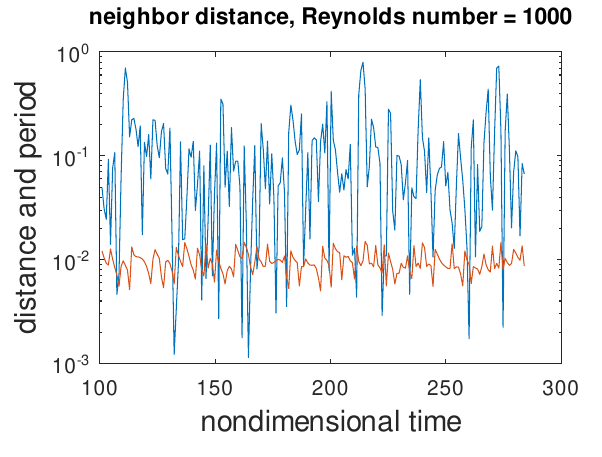}}
\vspace{-2mm}
\caption{Near neighbor distance $L$ for various Reynolds numbers, in blue.
The red curve is the computed approximate period $p(t)=t'-t$
multiplied by $10^{-3}$ for reference. Simulation data is the same as Figure \ref{fig:relfdata}.
}
\label{fig:nearnbraos1} 
\end{figure}

{The next parameter that we need to define is the distance between orbits.} We define this to be $L(t_0)=|x(t_0) - x(t_0')|$, where the vertical bars indicate
Euclidean distance.
In Figure \ref{fig:nearnbraos1}, we plot $L(t)$ as a function of $t$
as well as the aproximate period $p(t)=t'-t$ for various Reynolds numbers {(note that $p(t)$ is equivalent to one of the values $p_i$ in the iterative process described above, with $\sigma$ taken to be very small)}. For, Reynolds number 200, we see that $L$ is quite small
and decreasing as $t$ increases, as we would expect as the phase diagram
approaches a periodic orbit.
As the Reynolds number increases, $L$ increases in size, and for $R\geq 450$,
it no longer decreases as $t$ increases, indicating that the phase diagram
is no longer close to periodic.
We have also included in Figure \ref{fig:nearnbraos} plots of the computed
period $p(t)=t'-t$, scaled by $10^{-3}$ to fit on the $L$ plot.
As the Reynolds number increases, $p$ begins to oscillate, again indicating
a departure from periodicity.

The failure of the algorithm to identify the period
does not affect the computation of the Lyapunov exponent since
all that is required is to have a data point near the original one.
This means that we want $|L(t_i)|$ to be small.
As indicated in Figure \ref{fig:nearnbraos1}, this typically happens even 
if the point displacement differs from the average period substantially.
We see this by looking at the average distances (lower plot) in 
Figure \ref{fig:nearnbraos1}, which are typically less that 0.1.

\subsection{The parameter $\Delta$ and the approximation of $\lambda_1$}

The third and final parameter for the algorithm  \cite[section 5.1]{ref:timeseriesLyapunov}
is a small time parameter $\Delta$, not to be confused with the
time step $\Delta t$.
The notational confusion is significant, but we kept the notation $\Delta$
to remain as close to \cite{ref:timeseriesLyapunov} as possible. {Note that $\Delta$ operates in the same way as $\sigma$ in the computation of $p_S$}

{With $\Delta$ defined, we may finally describe the method for approximating $\lambda_1$.}
Consider the iterative nearest neighbor iteration process described in Appendix \ref{avP}, where the approximate period used is the best approximation, $p_S$, and the time delay between samples is $\Delta$.  

{For this process, we have $t_{k}=t_0+k\Delta=t_{k-1}+\Delta$, and at each iteration we may define}
$$
L'(t_{k})=|x(t_k)-x(t_{k-1}'+\Delta)|,\qquad
L(t_{k})=|x(t_k)-x(t_k')|,
$$
where again $t_k'$ is determined so that $x(t_k)$ and $x(t_k')$ are nearest neighbors.

The largest Lyapunov exponent $\lambda_1$ is then approximated as
\begin{equation}\label{lyapSum}
\lambda_1=\frac{1}{t_M-t_0}\sum_{k=1}^N \log_2 \frac{L'(t_{k})}{L(t_{k-1})},
\end{equation}
where $N$ is the number of $t_k$ values taken to approximate $\lambda_1$. We note that this number $N$ will depend on the size of $\Delta$.

{Note that we can write $L'(t_k) = \vert x(t_{k-1} + \Delta) - x(t_{k-1'}' + \Delta)\vert$. Consequently, the quotient in \eqref{lyapSum} may be written}
\begin{equation*}
\frac{L'(t_{k})}{L(t_{k-1})} = \frac{\vert x(t_{k-1} + \Delta) - x(t_{k-1'}' + \Delta)\vert}{\vert x(t_{k-1}) - x(t_{k-1}')\vert}.
\end{equation*}
{Thus, this ratio measures how the nearest neighbor evolves in time. If $x(t_{k-1}')$ remains close to $x(t_{k-1})$ as they both evolve,which we expect from a nearly periodic system, then $L'(t_{k})$ should be close to $L(t_{k-1})$. In this case, the ratio will be close to $1$, so its logarithm will be near zero, keeping $\lambda_1$ small.  Similarly, if the system is in the process of becoming more periodic/stabilizing, then we expect the trajectories to become closer together,  making $L(t_{k-1})$ the larger term in the ratio. This will result in a negative logarithm, and a more negative value of $\lambda_1$.}

{On the other hand, if the system is more chaotic, we expect the trajectories of $x(t_{k-1})$ and $x(t_{k-1}')$ to diverge. In particular, this will mean $L'(t_{k})$ is the larger term, so the logarithm will be positive, making $\lambda_1$ more positive. }

It can happen that $L'(t_k)=0$ or $L(t_{k-1})=0$, and $\lambda_1$ is undefined.
In our implementation, we define the quotient to be 1 in either case, essentially
skipping that time point.
However, such exceptions are rare, especially for $m>1$.

The effect of the parameters $m$ and $\Delta$ can be seen in 
Table \ref{tabl:mleo1}.
For Reynolds numbers $R=60, 120$ (and other not shown), we see that it is important to choose $m>1$,
but for $m\geq 2$ there is little change.
Thus we took $m=5$ for typical computations.
We do not fully understand why it is important to take $m>1$.

The dependence on $\Delta$ is more complicated, as shown in
Table \ref{tabl:mleo1} for Reynolds numbers $R=60$ and  $R=120$.
For the vortex street, the Strouhal period is roughly 10, so $\Delta =1$
means we are interrogating roughly 10\% of a period.
There is a basic trade-off regarding the choice of $\Delta$.
For larger $\Delta$, there are fewer intervals
available to define $\lambda_1$.
Correspondingly, for smaller $\Delta$, the
interrogation interval is smaller and thus less reliable.
We chose $\Delta =1$ as the best compromise. {In Table \ref{tabl:mleo1}, we indicate the dependence of the resulting estimate
of the Lyapunov exponent $\lambda_1$ on the various parameters of the model. }

\begin{table}
\begin{center} 
\begin{tabular}{|c|c|c|c|c|c|}
\hline
$R$ & $\lambda_1$ & $\Delta$ &start &end&m\\
\hline
           60 &    0.74161   &      0.2     &    380    &     480 &          1\\
           60 &    0.11183   &      0.2     &    380    &     480 &          2\\
           60 &    0.10191   &      0.2     &    380    &     480 &          3\\
           60 &    0.33788   &      0.5     &    380    &     480 &          1\\
           60 &    0.07838   &      0.5     &    380    &     480 &          2\\
           60 &    0.080907  &      0.5     &    380    &     480 &          3\\
           60 &    0.072739  &      0.5     &    380    &     480 &          4\\
           60 &    0.22372   &        1     &    380    &     480 &          1\\
           60 &    0.024088  &        1     &    380    &     480 &          2 \\ 
           60 &    0.025619  &        1     &    380    &     480 &          3 \\ 
           60 &    0.025915  &        1     &    380    &     480 &          5\\
           60 &    0.11994   &        2     &    380    &     480 &          1\\
           60 &    0.0016816 &        2     &    380    &     480 &          2\\
           60 &    0.0016274 &        2     &    380    &     480 &          3\\
           60 &    0.019166  &        3     &    380    &     480 &          3\\
\hline
\end{tabular}
\begin{tabular}{|c|c|c|c|c|c|}
\hline
$R$ & $\lambda_1$ & $\Delta$ &start &end&m\\
\hline
          120 &    0.14822   &      0.1     &    380    &     480 &1\\
          120 &    0.16902   &      0.2     &    380    &     480 &1\\
          120 & -0.0018047   &      0.5     &    380    &     480 &          1\\
          120 & -0.0018571   &      0.5     &    380    &     480 &          2\\
          120 & -0.0018931   &      0.5     &    380    &     480 &          3\\
          120 & -0.0019278   &      0.5     &    380    &     480 &          4\\
          120 &-0.00037597   &        1     &    380    &     480 &          1\\
          120 &-0.00043043   &        1     &    380    &     480 &          2\\
          120 &-0.00044859   &        1     &    380    &     480 &          3\\
          120 &-0.00045767   &        1     &    380    &     480 &          4\\
          120 &-0.00047329   &        1     &    380    &     480 &          5\\
          120 & -0.0010856   &        2     &    380    &     480 &1\\
          120 & -0.0011026   &        2     &    380    &     480 &          3\\
          120 &   0.042537   &      0.5     &    280    &     480 &1\\
          120 &   0.039571   &      0.2     &    280    &     480 &1\\
\hline
\end{tabular}
\end{center}
\vspace{-6mm}
\caption{Maximal Lyapunov exponent $\lambda_1$ as computed by the algorithm 
in Section \ref{sec:imexeme}.
$R$ is the Reynolds number, $m$ is the delay coordinate index, start and end
indicate the simulation time interval analyzed,
and $\Delta$ is the evolution time.
Meshsize $M=32$,  segments = 2048, time step $\Delta t=0.1$.
}
\label{tabl:mleo1}
\end{table}


\begin{thebibliography}{10}

\bibitem{alnaes2015fenics}
Martin Aln{\ae}s, Jan Blechta, Johan Hake, August Johansson, Benjamin Kehlet,
  Anders Logg, Chris Richardson, Johannes Ring, Marie~E. Rognes, and Garth~N.
  Wells.
\newblock The {FEniCS} project version 1.5.
\newblock {\em Archive of Numerical Software}, 3(100), 2015.

\bibitem{ref:BeaudanMointhesis94}
Patrick~Bruno Beaudan.
\newblock {\em Numerical experiments on the flow past a circular cylinder at
  sub-critical {Reynolds} number}.
\newblock PhD thesis, Stanford University, 1994.

\bibitem{birkhoff1953formation}
Garrett Birkhoff.
\newblock Formation of vortex streets.
\newblock {\em Journal of Applied Physics}, 24(1):98--103, 1953.

\bibitem{chen1995bifurcation}
J.-H. Chen, W.~G. Pritchard, and S.~J. Tavener.
\newblock Bifurcation for flow past a cylinder between parallel planes.
\newblock {\em Journal of Fluid Mechanics}, 284:23--41, 1995.

\bibitem{ref:DelaneySorensenNACAcyldrag}
Noel~K. Delany and Norman~E. Sorensen.
\newblock Low-speed drag of cylinders of various shapes.
\newblock Technical report, National Advisory Committee for Aeronautics, 1953.

\bibitem{ref:dongkarni05dns}
Suchuan Dong and George~E. Karniadakis.
\newblock {DNS of flow past a stationary and oscillating cylinder at Re}=
  10000.
\newblock {\em Journal of fluids and structures}, 20(4):519--531, 2005.

\bibitem{ref:dryden1930pressure}
Hugh~L. Dryden and George~C. Hill.
\newblock The pressure of the wind on large chimneys.
\newblock {\em Proceedings of the National Academy of Sciences},
  16(11):727--731, 1930.

\bibitem{ref:wolfchaocylinderflow}
D.~Durante, C.~Pilloton, and A.~Colagrossi.
\newblock Intermittency patterns in the chaotic transition of the planar flow
  past a circular cylinder.
\newblock {\em Physical Review Fluids}, 7(5):054701, 2022.

\bibitem{ref:dragecrisisfage}
Arthur Fage.
\newblock {XXVIII. The} air-flow around a circular cylinder in the region where
  the boundary layer separates from the surface.
\newblock {\em The London, Edinburgh, and Dublin Philosophical Magazine and
  Journal of Science}, 7(42):253--273, 1929.

\bibitem{ref:FeyStrouhalFormula}
Uwe Fey, Michael K{\"o}nig, and Helmut Eckelmann.
\newblock A new {Strouhal--Reynolds-number} relationship for the circular
  cylinder in the range {$47<$} {Re}{$< 2\times 10^5$}.
\newblock {\em Physics of Fluids}, 10(7):1547--1549, 1998.

\bibitem{ref:LinkeBeltramiFlows}
Nicolas~R. Gauger, Alexander Linke, and Philipp~W. Schroeder.
\newblock On high-order pressure-robust space discretisations, their advantages
  for incompressible high {Reynolds number generalised Beltrami} flows and
  beyond.
\newblock {\em arXiv preprint arXiv:1808.10711}, 2018.

\bibitem{lrsBIBiw}
Ingeborg~G. Gjerde and L.~Ridgway Scott.
\newblock Kinetic-energy instability of flows with slip boundary conditions.
\newblock {\em Journal of Mathematical Fluid Dynamics}, 24:97, 2022.

\bibitem{lrsBIBkj}
Ingeborg~G. Gjerde and L.~Ridgway Scott.
\newblock Evaluation of drag in fluid flow.
\newblock {\em submitted}, 2023.

\bibitem{lrsBIBjn}
Ingeborg~G. Gjerde and L.~Ridgway Scott.
\newblock Resolution of {D'Alembert's Paradox} using {Navier's} slip boundary
  conditions.
\newblock submitted, 2023.

\bibitem{ref:FazleHayakawatreedee}
Michio Hayakawa and Fazle Hussain.
\newblock Three-dimensionality of organized structures in a plane turbulent
  wake.
\newblock {\em Journal of Fluid Mechanics}, 206:375--404, 1989.

\bibitem{ref:cylendardrag}
C.~F. Heddleson, D.~L. Brown, and R.~T. Cliffe.
\newblock Summary of drag coefficients of various shaped cylinders.
\newblock Technical report, General Electric Co., Cincinnati OH, 1957.

\bibitem{ref:cylinderHopfvortexstreet}
C.~P. Jackson.
\newblock A finite-element study of the onset of vortex shedding in flow past
  variously shaped bodies.
\newblock {\em Journal of Fluid Mechanics}, 182:23--45, 1987.

\bibitem{ref:confinedcylinderflow}
Nicolas Kanaris, Dimokratis Grigoriadis, and Stavros Kassinos.
\newblock Three dimensional flow around a circular cylinder confined in a plane
  channel.
\newblock {\em Physics of Fluids}, 23(6):064106, 2011.

\bibitem{keulegan1956forces}
Garbis~H. Keulegan and Lloyd~H. Carpenter.
\newblock {\em Forces on cylinders and plates in an oscillating fluid},
  volume~60.
\newblock National Bureau of Standards, 1958.

\bibitem{ref:fracordimeapp}
Krishna Kumaraswamy.
\newblock Fractal dimension for data mining.
\newblock {\em Center for Automated Learning and Discovery School of Computer
  Science Carnegie Mellon University}, 5000, 2003.

\bibitem{ref:midRecylinderdrag}
John~H. Lienhard.
\newblock Synopsis of lift, drag, and vortex frequency data for rigid circular
  cylinders.
\newblock Technical report, 1966.

\bibitem{ref:Meneghini1994Thesis}
Julio~Romano Meneghini.
\newblock {\em Numerical simulation of bluff body flow control using a discrete
  vortex method}.
\newblock PhD thesis, Imperial College London, 1993.

\bibitem{lrsBIBia}
Hannah Morgan and L.~Ridgway Scott.
\newblock Towards a unified finite element method for the {Stokes} equations.
\newblock {\em SIAM Journal on Scientific Computing}, 40(1):A130--A141, 2018.

\bibitem{oliver2001vortex}
Marcel Oliver and Steve Shkoller.
\newblock The vortex blob method as a second-grade {non-Newtonian} fluid.
\newblock {\em Communications in Partial Differential Equations}, 26:295–314,
  2001.

\bibitem{PantonRonaldL2013IF}
Ronald~L. Panton.
\newblock {\em Incompressible Flow}.
\newblock John Wiley {\&} Sons, Incorporated, Somerset, fourth edition, 2013.

\bibitem{ref:Rajani2D3Dcylinderflow}
B.~N. Rajani, A.~Kandasamy, and Sekhar Majumdar.
\newblock Numerical simulation of laminar flow past a circular cylinder.
\newblock {\em Applied Mathematical Modelling}, 33(3):1228--1247, 2009.

\bibitem{ref:errelflowRex}
E.~R. Relf.
\newblock Discussion of the results of measurements of the resistance of wires,
  with some additional tests on the resistance of wires of small diameter.
\newblock Technical Report 102, Reports and Memoranda of the Aeronautical
  Research Council, 1914.

\bibitem{lrsBIBih}
L.~Ridgway Scott.
\newblock {\em Introduction to Automated Modeling with FEniCS}.
\newblock Computational Modeling Initiative, 2018.

\bibitem{lrsBIBis}
L.~Ridgway Scott.
\newblock Kinetic energy flow instability with application to {Couette} flow.
\newblock Research Report {UC/CS} TR-2020-07, Dept. Comp. Sci., Univ. Chicago,
  2020.

\bibitem{ref:SkokosreviewArxiv08}
Charalampos Skokos.
\newblock The {Lyapunov} characteristic exponents and their computation.
\newblock In Jean~J. Souchay and Rudolf Dvorak, editors, {\em Dynamics of Small
  Solar System Bodies and Exoplanets}, pages 63--135. Springer, 2010.

\bibitem{ref:karmanstreetexpts}
Sadatoshi Taneda.
\newblock Experimental investigation of vortex streets.
\newblock {\em Journal of the Physical Society of Japan}, 20(9):1714--1721,
  1965.

\bibitem{ref:trittonlowrexpts}
David~J. Tritton.
\newblock Experiments on the flow past a circular cylinder at low {Reynolds}
  numbers.
\newblock {\em Journal of Fluid Mechanics}, 6(4):547--567, 1959.

\bibitem{ref:vanDyke}
Milton Van~Dyke and Milton Van~Dyke.
\newblock {\em An album of fluid motion}, volume 176.
\newblock Parabolic Press Stanford, 1982.

\bibitem{wieselsberger1921neuere}
C.~von Wieselsberger.
\newblock Neuere feststellungen uber die gesetze des flussigkeits und
  luftwiderstands.
\newblock {\em Phys. Z.}, 22:321, 1921.

\bibitem{ref:WieselsbergerNACA}
Carl Wieselsberger.
\newblock New data on the laws of fluid resistance.
\newblock Technical Report NACA Technical Note No.~84, 1922.

\bibitem{wiggins2003introduction}
Stephen Wiggins.
\newblock {\em Introduction to applied nonlinear dynamical systems and chaos}.
\newblock Springer, 2nd edition, 2003.

\bibitem{ref:annurevfluidvortexvibrations}
Charles H.~K. Williamson and R.~Govardhan.
\newblock Vortex-induced vibrations.
\newblock {\em Annu. Rev. Fluid Mech.}, 36:413--455, 2004.

\bibitem{ref:timeseriesLyapunov}
Alan Wolf, Jack~B. Swift, Harry~L. Swinney, and John~A. Vastano.
\newblock Determining {Lyapunov} exponents from a time series.
\newblock {\em Physica D: nonlinear phenomena}, 16(3):285--317, 1985.

\end{thebibliography}
\end{document}